\documentclass[a4paper, 11pt]{article}
\usepackage{amsmath,amssymb,amsthm,xypic} 
\input xy
\xyoption{all}

\title{A Chevalley theorem for difference equations}
\author{Michael Wibmer}

\newtheorem{theo}{Theorem}[section]
\newtheorem{lemma}[theo]{Lemma}
\newtheorem{prop}[theo]{Proposition}
\newtheorem{cor}[theo]{Corollary}
\newtheorem{defi}[theo]{Definition}
\newtheorem{rem}[theo]{Remark}

\theoremstyle{definition}
\newtheorem{ex}[theo]{Example}

\newcommand{\ida}{\mathfrak{a}}

\newcommand{\p}{\mathfrak{p}}
\newcommand{\q}{\mathfrak{q}}
\newcommand{\m}{\mathfrak{m}}

\newcommand{\spec}{\operatorname{Spec}}

\newcommand{\gal}{\operatorname{Gal}}
\newcommand{\Gl}{\operatorname{GL}}

\newcommand{\Q}{\mathfrak{Q}}

\newcommand{\id}{\operatorname{id}}

\newcommand{\s}{\sigma}
\newcommand{\de}{\delta}
\newcommand{\ssp}{\sigma^\bullet\text{-}\operatorname{Spec}}

\newcommand{\Spec}{\text{-}\operatorname{Spec}}

\begin{document}

\maketitle

\begin{abstract}
By a theorem of Chevalley the image of a morphism of varieties is a constructible set. The algebraic version of this fact is usually stated as a result on  ``extension of specializations'' or ``lifting of prime ideals''. We present a difference analog of this theorem. The approach is based on the philosophy that occasionally one needs to pass to higher powers of $\s$, where $\s$ is the endomorphism defining the difference structure. In other words, we consider difference pseudo fields (which are finite direct products of fields) rather than difference fields. We also prove a result on compatibility of pseudo fields and present some applications of the main theorem, e.g. constrained extension and uniqueness of differential Picard-Vessiot rings with a difference parameter.
\end{abstract}

\section*{Introduction}
A classical theorem of Chevalley states the following:
\begin{theo}
\label{theo usual geometric Chevalley theorem}
Let $f:X\rightarrow Y$ be a morphism of finite type between noetherian schemes. Then $f$ maps constructible sets to constructible sets. In particular the image $f(X)$ of $f$ is a constructible subset of $Y$.
\end{theo}

One can find several proofs of this theorem in the literature (e.g. \cite[Theorem 1.8.4, p. 239]{Grothendieck:EGAIV1}). Most of them rely on the following algebraic version of Chevalley's Theorem (which does not require the noetherianity assumption).

\begin{theo}[Algebraic Chevalley theorem]
 \label{theo usual algebraic chevalley}
Let $R\subset S$ be an inclusion of integral domains such that $S$ is finitely generated over $R$. Then there exists a non-zero element $r\in R$ such that every prime ideal $\q$ of $R$ with $r\notin \q$ lifts to a prime ideal of $S$, i.e. there exists a prime ideal $\q'$ of $S$ with $\q'\cap R=\q$.
\end{theo}
The geometric meaning of Theorem \ref{theo usual algebraic chevalley} is simply that the image of $f:\spec(S)\rightarrow\spec(R)$ contains the non-empty open subset $D(r)$.
The following version of Chevalley's Theorem is also quite popular. It is immediately seen to be equivalent to Theorem \ref{theo usual algebraic chevalley} via Hilbert's Nullstellensatz.

\begin{theo}
\label{theo usual algebraic chevalley with closed field}
Let $R\subset S$ be an inclusion of integral domains such that $S$ is finitely generated over $R$. Then there exists a non-zero element $r\in R$ such that every morphism $\psi:R\rightarrow k$ into an
algebraically closed field $k$ with $\psi(r)\neq 0$ can be extended to a morphism $\widetilde{\psi}: S\rightarrow k$.
\end{theo}

The differential analogs of the above three theorems are all true. This is essentially Theorem 3, p. 140 in Kolchin's book \cite{Kolchin:differentialalgebraandalgebraicgroups}, who gives credit to Ritt, Seidenberg and Rosenfeld. There has been some recent interest in the differential version of Theorem \ref{theo usual algebraic chevalley with closed field} in connection with differentially closed fields, see \cite{Kac:differentialChevalley} and \cite{Rosen:differentialChevalley}.

As in many cases, the situation is slightly more complicated if we consider difference equations. First of all, it is well known (see \cite[Example 3, p. 214]{Cohn:difference} or Example \ref{ex failure of naive chevalley} below) that the ``naive'' difference analog of Theorem \ref{theo usual algebraic chevalley} fails. It appears that the best known approximations to a difference algebraic Chevalley theorem are the following two theorems of Cohn \cite[Theorem 11, p. 227 and Theorem 12\footnote{The characteristic zero assumption seems to be a missprint.}, p. 230]{Cohn:difference}.

\begin{theo}[Cohn]
Let $k$ be a $\s$-field and $R\subset S$ an inclusion of finitely $\s$-generated $k$-$\s$-algebras such that $S$ is an integral domain and $\s$ is injective on $S$.
Let $K$ and $L$ denote the quotient fields of $R$ and $S$ respectively and let $b$ be a finite tuple of generators of the core $L_K\subset L$ of $L$ over $K$. Then there exists a non-zero $r\in R$ with the following property: Every $\s$-prime ideal $\q$ of $R$ with $r\notin\q$ that lifts to a $\s$-prime ideal of $R\{b\}$ also lifts to a $\s$-prime ideal of $S$.
\end{theo}

Here, by a $\s$-prime ideal $\q$ of a difference ring (or $\s$-ring for short) we mean a difference ideal that is prime and reflexive, i.e $\s^{-1}(\q)=\q$. As the core $L_K$ is a separable algebraic extension of $K$ this implies that almost every $\s$-prime ideal of $R$ lifts to $S$ if $K$ is relatively separably algebraically closed in $L$.

\begin{theo}[Cohn]
Let $k$ be a $\s$-field and $R$ a $\s$-polynomial ring in finitely many $\s$-indeterminates  over $k$. Let $S$ be a finitely $\s$-generated $k$-$\s$-algebra containing $R$ such that $S$ is an integral domain and $\s$ is injective on $S$. Then there exists a non-zero $r\in R$ with the following property: Every $\s$-prime ideal $\q$ of $R$ with $r\notin \q$ such that the quotient field of $S$ and the residue field at $\q$ are compatible $\s$-field extensions of $k$ has a lift to a $\s$-prime ideal of $S$.
\end{theo}

To motivate our difference version of Chevalley's Theorem we first look at the standard counterexample.

\begin{ex} \label{ex failure of naive chevalley}
We consider the ring $\mathbb{C}[x]$ of univariate polynomials over $\mathbb{C}$ as difference ring by setting $\s(x)=-x$. Then $\mathbb{C}[x^2]\subset \mathbb{C}[x]$ is an inclusion of $\s$-rings. The non-zero $\s$-prime ideals of $R=\mathbb{C}[x^2]$ are precisely those of the form $\q=(x^2-a)$ for some $a\in\mathbb{C}$. On the other hand the only non-zero $\s$-prime ideal of $\mathbb{C}[x]$ is $(x)$. We see that $(0)$ and $(x^2)$ are the only $\s$-prime ideals of $\mathbb{C}[x^2]$ that lift to $\mathbb{C}[x]$. The heart of the problem is that although the fibre ring over $\q=(x^2-a)$ for $a\in\mathbb{C}\smallsetminus\{0\}$ is non-zero it has empty difference spectrum. However, note that $\q=(x^2-a)$ has two lifts to a $\s^2$-prime ideal of $S=\mathbb{C}[x]$, namely $(x-\sqrt{a})$ and $(x+\sqrt{a})$ which are permuted by $\s$.
\end{ex}
Our main result asserts that one can always resolve the situation by allowing a passage to higher powers of $\s$ as in the above example. Moreover we are able to prove a result that is uniform in powers of $\s$.

\begin{theo}[Difference algebraic Chevalley theorem] \label{theo difference algebraic Chevalley Introduction}
Let $R\subset S$ be an inclusion of $\s$-rings such that $S$ is an integral domain and $\s$ is injective on $S$. Assume that $S$ is finitely $\s$-generated over $R$. Then there exists an integer
$l\geq 1$ and a non-zero element $r\in R$ such that for all $d\geq 1$ and all $\s^d$-prime ideals $\q$ of $R$ with $r\s(r)\cdots\s^{d-1}(r)\notin\q$ there exists a $\s^{ld}$-prime ideal $\q'$ of $S$ with $\q'\cap R=\q$. In particular every $\s$-prime ideal $\q$ of $R$ with $r\notin\q$ lifts to a $\s^l$-prime ideal of $S$.
\end{theo}

As seen in the above example it is in general not possible to choose $l=1$. This is a special ``feature'' of difference algebra that appears to have no differential analog. However, we note that the philosophy of considering higher powers of a derivation was the key idea to generalize the main results of differential algebra to positive characteristic (see e.g. \cite{Okugawa:DifferentialAlgerbaofnonzeroCharacteristic} or \cite{Matzat:skript}). The idea of considering higher powers of $\s$ has already been successfully applied in \cite{SingerPut:difference}, \cite{Hrushovski:elementarytheoryoffrobenius} and \cite{Hrushovskietal:ModelTheoryofDifferenceFieldsIIPeriodicIdelas}. However the approach of \cite{Hrushovskietal:ModelTheoryofDifferenceFieldsIIPeriodicIdelas} appears to be somehow dual to the one presented here.

Note that by taking $\s$ equal to the identity in Theorem \ref{theo difference algebraic Chevalley Introduction} we precisely recover Theorem \ref{theo usual algebraic chevalley}.
Probably the most natural formulation of Theorem \ref{theo difference algebraic Chevalley Introduction} is in terms of $\s$-pseudo prime ideals and $\s$-pseudo fields. A $\s$-pseudo prime ideal $\p$ of a $\s$-ring $R$ is simply a difference ideal of the form $\p=\q\cap\s^{-1}(\q)\cap\cdots\cap\s^{-(d-1)}(\q)$ for some $\s^d$-prime ideal $\q$ of $R$. The total quotient ring of $R/\p$ is then a $\s$-pseudo field (See Section 1.1 for a more methodical definition).

\begin{cor} \label{cor chevalley with pseudo field introduction}
Let $R\subset S$ be an inclusion of $\s$-rings such that $S$ is an integral domain and $\s$ is injective on $S$. Assume that $S$ is finitely $\s$-generated over $R$.
Then there exists a non-zero $r\in R$ with the following property: Every morphism $\psi:R\rightarrow K$ into some $\s$-pseudo field $K$ such that $\psi(r)$ is invertible can be extended to a morphism $\widetilde{\psi}:S\rightarrow\widetilde{K}$ where $\widetilde{K}$ is a $\s$-pseudo field extension of $K$.
\end{cor}

Akin to Theorem \ref{theo usual algebraic chevalley with closed field} we can also give a result with $\s$-closed fields, i.e models of ACFA (See \cite{Hrushovskietal:ModelTheoryofDifferencefields}).

\begin{cor} \label{cor chevalley with closed field introduction}
Let $R\subset S$ be an inclusion of $\s$-rings such that $S$ is an integral domain and $\s$ is injective on $S$. Assume that $S$ is finitely $\s$-generated over $R$.
Then there exists a non-zero $r\in R$ with the following property: For every $\s$-morphism $\psi:R\rightarrow k$ into a $\s$-closed field $k$ with $\psi(r)\neq 0$ there exists an integer $l\geq1$ and an extension $\widetilde{\psi}:S\rightarrow k$ of $\psi$ to a $\s^l$-morphism on $S$.
\end{cor}

Results similar to Corollary \ref{cor chevalley with pseudo field introduction} and \ref{cor chevalley with closed field introduction} above have also been obtained
in \cite[Section 4.4]{Trushin:DifferenceNullstellsatzCaseOfFiniteGroup} when the action of $\s$ is replaced by the action of a finite group. But note that under this assumption finitely difference generated implies finitely generated and so the usual Chevalley theorem is available.

\vspace{5mm}

We also present some applications of Theorem \ref{theo difference algebraic Chevalley Introduction}. Most prominently we show that the Picard-Vessiot theory of a linear differential equation with a difference parameter is a perfect match with the philosophy of passing to higher powers of $\s$. For example, we prove that two Picard-Vessiot rings for the same differential equation over a $\de\s$-field $k$ with $\s$-closed $\de$-constants are isomorphic as $k$-$\de\s^l$-algebras for some integer $l\geq 1$. We also provide an example to illustrate that one can not choose $l=1$ in general. It is the hope of the author that this result will convince the reader, if he is still reluctant to accept pseudo-fields instead of fields, that it is worthwhile to adopt the new point of view.

\vspace{5mm}

The article is divided into two parts. The proof of Theorem \ref{theo difference algebraic Chevalley Introduction} is contained in the first part. In the second part we present the applications. The first part starts with fixing the notation and basic properties of $\s$-pseudo fields. For example we prove that any two extensions of the same $\s$-pseudo field are compatible if one of them is finitely $\s$-generated. Then, after some preparatory results, the proof of Theorem \ref{theo difference algebraic Chevalley Introduction} begins. The first part finishes with some variations of the main theorem and an explanation why we think that the result is more or less optimal.

The second part starts with a section on constrained extensions, whose theory heavily relies on Theorem \ref{theo difference algebraic Chevalley Introduction}. Philosophically constrained extensions play the same role in difference algebraic geometry as algebraic extensions in algebraic geometry.
The last section is on the Galois theory of linear differential equations with a difference parameter. We do not give an exposition of the whole theory -- this will be done in \cite{HardouindiVizio:DifferenceGaloisofDifferential}. The main result presented here says that two Picard-Vessiot rings with a $\s$-parameter for the same differential equation become isomorphic over a finitely $\s$-generated constrained extension of the differential constants.

\section{Pseudo fields and Chevalley theorem}

\subsection{Notation and preliminaries}

All rings are assumed to be commutative with identity. An ideal is called proper if it is not the whole ring or the zero ideal. For a ring $R$ we let $\Q(R)$ denote the total ring of quotients of $R$, i.e. the localization of $R$ at the multiplicatively closed subset of all non-zero divisors. For a prime ideal $\q$ of $R$ we denote with $k(\q)=\Q(R/\q)$ the residue field at $\q$.

A difference ring \footnote{We do not assume that $\s$ is injective.} (or $\s$-ring for short) is a ring $R$ together with an endomorphism $\s:R\rightarrow R$. Throughout the prefix ``$\s\text{-}$'' is to be read as ``difference'' or ``transformal''. We largely use standard notations of difference algebra as they can be found in \cite{Cohn:difference} and \cite{Levin:difference}. Following \cite{Hrushovski:elementarytheoryoffrobenius} we call a prime ideal $\q$ of $R$ with $\s^{-1}(\q)=\q$ a \emph{$\s$-prime ideal} \footnote{This is not in accordance with \cite[Definition 2.3.20, p. 128]{Levin:difference} or \cite[Definition 1.4.3, p. 11]{Wibmer:thesis}}. A \emph{$\s$-domain} is a $\s$-ring whose zero ideal is $\s$-prime, i.e. an integral domain with $\s$ injective.
The set of all $\s$-prime ideals of $R$ is denoted with $\s\Spec(R)$.
We put $\s^0=\id$ and for a subset $S$ of $R$ and an integer $d\geq 1$ we set $\s^{-d}(S)=\{r\in R|\ \s^d(r)\in S\}$. We say that $R$ is a $\s^d$-ring to indicate that we consider the difference ring $(R,\s^d)$. E.g. a $\s^d$-prime ideal $\q$ of $R$ is a prime ideal with $\s^{-d}(\q)=\q$. We say that a prime ideal $\q$ of $R$ is a $\s^\bullet$-prime ideal if it is a $\s^d$-prime ideal for some $d\geq 1$. The set of all $\s^\bullet$-prime ideals of $R$ is denoted with $\s^\bullet\Spec(R)$.
If $\q$ is a $\s^\bullet$-prime ideal, the smallest integer $d\geq 1$ such that $\q$ is a $\s^d$-prime ideal is called the \emph{period} of $\q$. (It is the period of the orbit of $\q$ under the action of $\s$ on $\spec(R)$.)

A $\s$-ring $K$ is called a \emph{$\s$-pseudo field} if it is noetherian, $\s$-simple (i.e. $K$ has no proper $\s$-ideal) and every non-zero divisor of $K$ is invertible. As a ring with no proper ideals is a field one my regard $\s$-pseudo fields as the difference analog of fields. Pseudo fields have essentially been introduced by Singer and van der Put in \cite{SingerPut:difference} to avoid certain pathologies in the study of solutions of linear difference equations that arise if one restricts the attention to fields. One aim of this article is to show that the underlying idea of pseudo fields can also be applied to the study of general algebraic difference equations. Pseudo fields have also been used in \cite{AmanoMasuoka:artiniansimple}, \cite{Wibmer:thesis}, \cite{Trushin:DifferenceNullstellensatz} and \cite{Trushin:DifferenceNullstellsatzCaseOfFiniteGroup}. The notion of pseudo field in \cite{Trushin:DifferenceNullstellensatz} and \cite{Trushin:DifferenceNullstellsatzCaseOfFiniteGroup} is however more general.

If $K$ is a $\s$-pseudo field there exist uniquely determined idempotent elements $e_1,\ldots,e_d\in K$ such that
\begin{itemize}
 \item $e_ie_j=0$ for $i\neq j$ and $\s(e_i)=e_{i+1}$ for $i=1,\ldots,d$ (with $e_{d+1}:=e_1$),
\item $e_iK$ is a field (with identity element $e_i$) and
\item $K=e_1K\oplus\cdots\oplus e_dK$.
\end{itemize}
(See \cite[Corollary 1.16, p. 12]{SingerPut:difference} or \cite[Proposition 1.3.2, p. 9]{Wibmer:thesis}.)
We call $d$ the \emph{period of $K$}. Note that $e_iK$ is a $\s^d$-field. If we say something like ``Let $K=e_1K\oplus\cdots\oplus e_dK$ be a $\s$-pseudo field.'' we always mean that the $e_i$'s are as specified above.

One easily shows (\cite[Proposition 1.4.2]{Wibmer:thesis}) that for a $\s$-ideal $\p$ of a $\s$-ring $R$ the following statements are equivalent.
\begin{enumerate}
 \item $\p$ is the kernel of a $\s$-morphism $R\rightarrow K$ into a $\s$-pseudo field $K$.
\item $\s$ extends to $\Q(R/\p)$ and $\Q(R/\p)$ is a $\s$-pseudo field.
\item There exists a $\s^d$-prime ideal $\q$ of $R$ with
\[\p=\q\cap\s^{-1}(\q)\cap\cdots\cap \s^{-(d-1)}(\q).\]
\end{enumerate}

A $\s$-ideal satisfying the above properties will be called a \emph{$\s$-pseudo prime ideal}. In \cite{Trushin:DifferenceNullstellensatz} one can find a more general notion of pseudo prime ideal. The $\s$-pseudo field $k(\p)=\Q(R/\p)$ is called the residue $\s$-pseudo field at $\p$.

We call a $\s$-ring a \emph{$\s$-pseudo domain} if the zero ideal is $\s$-pseudo prime. Any $\s$-subring of a $\s$-pseudo domain is a $\s$-pseudo domain and for any $\s$-subring $R$ of a $\s$-pseudo field $K$ the $\s$-pseudo field $\Q(R)$ is naturally embedded into $K$. Let $L|K$ be an extension of $\s$-pseudo fields and $S$ a subset of $L$. We let
$K\langle S\rangle=\Q(K\{S\})\subset L$ denote the smallest $\s$-pseudo subfield of $L$ that contains $K$ and $S$. We say that $L$ is finitely $\s$-generated over $K$ (as $\s$-pseudo field) if there exists a finite subset $S$ of $L$ such that $L=K\langle S\rangle$.

Let $R$ be a $\s$-ring and $r\in R$. Let $T$ denote the multiplicatively closed subset of $R$ generated by $r,\s(r),\ldots$. We denote the localization $T^{-1}R$ with $R\{\frac{1}{r}\}$ (It is naturally a $\s$-ring). For a $\s$-ideal $\ida$ of $R$ we define $r\notin_\s\ida$ to mean $T\cap\ida=\emptyset$ and $\in_\s$ to be the negation of $\notin_\s$. Note that for a $\s$-prime ideal $\q$ we have $r\notin_\s\q$ if and only if $r\notin\q$.

\subsection{The compatibility theorem}

One of the disagreeable phenomenons in difference algebra is that even quite simple difference rings may have empty difference spectrum. One of the aims of this article is to show that this can be partly compensated by considering $\s^\bullet\Spec(R)$ instead of $\s\Spec(R)$.
This has a simple interpretation in terms of equations: Let $k$ be a $\s$-field and let $k\{x\}=k\{x\}_\s$ denote the $\s$-polynomial ring in the $\s$-variables $x=(x_1,\ldots,x_n)$ over $k$. Then, for any $d\geq 1$,
$k\{x\}_\s=k\{x,\s(x),\ldots,\s^{d-1}(x)\}_{\s^d}$ can equally well be interpreted as the $\s^d$-polynomial ring in the $\s^d$-variables $(x,\s(x),\ldots,\s^{d-1}(x))$ over the $\s^d$-field $k$. Thus a system of algebraic $\s$-equations $S\subset k\{x\}$ over the $\s$-field $k$ gives rise to a sequence $S_d$ of systems of algebraic $\s^d$-equations over the $\s^d$-field $k$. A $\s^d$-prime ideal of $k\{x\}$ containing $S$ corresponds to a solution of $S_d$.
The idea is that for understanding $S$ it is helpful to understand the whole family $(S_d)_{d\geq 1}$. Note that it can happen that $S$ has no solution (in a $\s$-overfield of $k$) but $S_d$ has a solution (in some $\s^d$-overfield of $k$).

It appears that almost all constructions of difference algebra are compatible with the philosophy of passing to higher powers of $\s$. The following example of this phenomenon will be needed below.

\begin{lemma} \label{lemma limit degree}
Let $L|K$ be a finitely $\s$-generated extension of $\s$-fields and $d\geq 1$ an integer. Then $L|K$ is finitely $\s^d$-generated (as an extension of $\s^d$-fields) and the limit degrees are related by the formula
\[\s^d\text{-}\operatorname{ld}(L|K)=(\s\text{-}\operatorname{ld}(L|K))^d.\]
Moreover the $\s$-core $L_K$ of $L|K$ equals the $\s^d$-core of $L|K$.
\end{lemma}
\noindent Proof: If $a=(a_1,\ldots,a_n)$ is a $\s$-generating set for $L|K$ then $b=(a,\s(a),\ldots,\s^{d-1}(a))$ is a $\s^d$-generating set. For $m\gg 1$ we have
\begin{align*}
\s^d\text{-}\operatorname{ld}(L|K)& =\left[K\big(b,\s^d(b),\ldots,\s^{md}(b)\big) : K\big(b,\ldots,\s^{(m-1)d}(b)\big)\right] \\
 & = \left[K\big(a,\s(a),\ldots,\s^{(m+1)d-1}(a)\big) : K\big(a,\ldots,\s^{md-1}(a\big)\right] \\
& =(\s\text{-}\operatorname{ld}(L|K))^d.
\end{align*}

Recall (e.g. \cite[Chapter 7, Section 15, p. 215]{Cohn:difference}) that the $\s$-core $L_K$ consists of all elements $a \in L$ that are separable algebraic over $K$ and satisfy $\s\text{-}\operatorname{ld}(K\langle a\rangle_\s|K)=1$. If $\s\text{-}\operatorname{ld}(K\langle a\rangle_\s|K)=1$ then also $\s^d\text{-}\operatorname{ld}(K\langle a\rangle_\s|K)=1$ and by the composition theorem for $\s^d\text{-}\operatorname{ld}$ we conclude that $\s^d\text{-}\operatorname{ld}(K\langle a\rangle_{\s^d}|K)=1$. Conversely if
$\s^d\text{-}\operatorname{ld}(K\langle a\rangle_{\s^d}|K)=1$ then also $\s^d\text{-}\operatorname{ld}(K\langle\s^i(a)\rangle_{\s^d}|K)=1$ for $i=0,\ldots,d-1$. As $K\langle a\rangle_\s$ is the compositum of $K\langle a\rangle_{\s^d},\ldots,K\langle\s^{d-1}(a)\rangle_{\s^d}$ it follows $\s^d\text{-}\operatorname{ld}(K\langle a\rangle_\s|K)=1$ and therefore $\s\text{-}\operatorname{ld}(K\langle a\rangle_\s|K)=1$. \qed

\vspace{5mm}

The problem of empty $\s$-spectrum is closely related to the classical problem of incompatibility (see \cite[Chapter 7]{Cohn:difference} or \cite[Chapter 5]{Levin:difference}). Recall that two extensions $L|K$ and $L'|K$ of $\s$-fields are said to be compatible (as $\s$-fields) if there exists a $\s$-field extension $M$ of $K$ containing both $L$ and $L'$ (up to $K$-$\s$-isomorphisms). Obviously $L|K$ and $L'|K$ are compatible if and only if $\s\Spec(L\otimes_K L')$ is not empty.

Following the above nomenclature we say that two $\s$-pseudo field extensions $L|K$ and $L'|K$ are \emph{compatible (as $\s$-pseudo fields)} if there exists a $\s$-pseudo field extension $M$ of $K$ with $K$-$\s$-embeddings $L\rightarrow M$, $L'\rightarrow M$. Now $L|K$ and $L'|K$ are compatible if and only if $\s^\bullet\Spec(L\otimes_K L')$ is not empty.

There is a certain class of hoped-for results in difference algebra that are not quite true, i.e. they are only available under certain additional assumptions. The basic idea is that these statements become valid unconditionally if one replaces difference fields with difference pseudo fields. To illustrate this idea we shall prove the following.

\begin{theo} \label{theo compatibility}
Let $K$ be a $\s$-pseudo field and $L,L'$ $\s$-pseudo field extensions of $K$. Assume that one of them is finitely $\s$-generated over $K$. Then $L|K$ and $L'|K$ are compatible.
\end{theo}
\noindent Proof: The proof is essentially a reduction to the classical compatibility theorem (\cite[Theorem 8, p. 223]{Cohn:difference}). We assume that $L$ is finitely $\s$-generated over $K$. As the first step we show that we can assume without loss of generality that $K$ is a field. Let $K=e_1K\oplus\cdots\oplus e_dK$. Then
\[L\otimes_K L'=(e_1L\otimes_{e_1K} e_1L')\oplus\cdots\oplus(e_dL\otimes_{e_dK} e_dL').\]
Note that $e_iL$ is a $\s^d$-pseudo field extension of $e_iK$ and $e_iL$ is finitely $\s^d$-generated over $e_iK$. As $\s$-permutes the idempotents $e_i\otimes e_i$ in a cyclic way we see that if $(e_1L\otimes_{e_1K} e_1L')$ has a $\s^{\bullet d}$-prime ideal then $L\otimes_K L'$ has a $\s^\bullet$-prime ideal. Therefore, we can assume that $K$ is a field.

Next we will reduce to the case that also $L$ is a field. If $L=e_1L\oplus\cdots\oplus e_dL$ then
\[L\otimes_K L'=(e_1L\otimes_K L')\oplus\cdots\oplus(e_dL\otimes_KL').\]
As $\s$-permutes the idempotents $e_i\otimes 1$ in a cyclic way we see that if $e_1L\otimes_K L'$ has a $\s^{\bullet d}$-prime ideal then $L\otimes_K L'$ has a $\s^\bullet$-prime ideal. Moreover $e_iL$ is finitely $\s^d$-generated over $K$. Consequently we can assume without loss of generality the $L$ is a field. In a similar fashion we can reduce to the case that also $L'$ is a field.

Thus we can assume from now that $L|K$ and $L'|K$ are $\s$-field extensions. We have to show that there is an integer $d\geq 1$ such that $L|K$ and $L'|K$ are compatible as $\s^d$-fields. As $L$ is finitely $\s$-generated over $K$ it follows (\cite[Theorem 18, p. 145 and Theorem 17, p. 144]{Cohn:difference}) that the core $L_K$ of $L|K$ is a finite algebraic extension of $K$. This implies that $L_K\otimes_K L'$ is noetherian. It is easy to see (cf. Remark \ref{rem pseudo simple} (i)) that a noetherian $\s$-ring has non-empty $\s^\bullet$-spectrum. It follows that $L_K\otimes_K L'$ has a $\s^d$-prime ideal for some $d\geq 1$, i.e. $L_K|K$ and $L'|K$ are compatible as $\s^d$-fields. By Lemma \ref{lemma limit degree} above the $\s$-core of $L|K$ equals the $\s^d$-core of $L|K$. It thus follows from the classical compatibility theorem \cite[Corollary, p.224]{Cohn:difference} that $L|K$ and $L'|K$ are compatible as $\s^d$-fields. \qed

\vspace{5mm}

For the reader reluctant to deal with pseudo fields we record the purely field theoretic version.
\begin{cor} \label{cor compatibility}
Let $L|K$ and $L'|K$ be $\s$-field extensions such that one of them is finitely $\s$-generated. Then there exists an integer $d\geq 1$ such that $L|K$ and $L'|K$ are compatible as $\s^d$-fields. \qed
\end{cor}

Unfortunately Theorem \ref{theo compatibility} is not true without the assumption of finite generation. This is illustrated in the following example.
\begin{ex}
Let $p$ be a prime number and let $K=\mathbb{F}_p$ denote the field with $p$ elements, considered as constant $\s$-field. Let $L=\overline{\mathbb{F}_p}$ denote the algebraic closure of $\mathbb{F}_p$ considered as difference field via the Frobenius, i.e. $\s(a)=\s_p(a)=a^p$ and let $L'=\overline{\mathbb{F}_p}$ denote the constant $\s$-field on $\overline{\mathbb{F}_p}$. Then $L\otimes_K L'$ has no $\s^\bullet$-prime ideal. This can be seen as follows:

The prime ideals of $L\otimes_K L'$ are all maximal (and minimal) and in bijection with the elements of the Galois group $\gal(\overline{\mathbb{F}_p}|\mathbb{F}_p)=\widehat{\langle \s_p\rangle}$.
If $\q$ is a prime ideal of $\overline{\mathbb{F}_p}\otimes_{\mathbb{F}_p}\overline{\mathbb{F}_p}$ then we have two isomorphisms
$\tau_s,\tau_t:\overline{\mathbb{F}_p}\rightarrow (\overline{\mathbb{F}_p}\otimes_{\mathbb{F}_p}\overline{\mathbb{F}_p})/\q$ defined by $\tau_s(a)=\overline{a\otimes 1}$ and
$\tau_t(a)=\overline{1\otimes a}$. Then $\tau_\q=\tau=\tau_s^{-1}\tau_t$ is an automorphism of $\overline{\mathbb{F}_p}$ over $\mathbb{F}_p$.
A $\s^d$-prime ideal $\q$ of $L\otimes_K L'$ would correspond to an automorphism such that $\tau_\q=\tau_\q\s_p^d$. But this is impossible as the Frobenius is not of finite order.
\end{ex}

The phenomenon of incompatibility can have the following awkward effect: A system of algebraic difference equations with coefficients in some $\s$-field $K$ that has a solution (in some $\s$-overfield of $K$) need not have a solution (in some $\s$-overfield of $K'$) if we think of the system as having coefficients in some $\s$-overfield $K'$ of $K$. In other words the base extensions of a non-empty $\s$-variety can be empty. The following geometric corollary to Theorem \ref{theo compatibility} is our proposed solution to this problem.

\begin{cor}
Let $K$ be a $\s$-field and $R$ a finitely $\s$-generated $K$-$\s$-algebra. Then for every $\s$-field extension $K'$ of $K$ the canonical map
\[ \ssp(R\otimes_K K')\longrightarrow \ssp(R)\]
is surjective. In particular $\ssp(R\otimes_K K')$ is non-empty if $\ssp(R)$ is non-empty.
\end{cor}
\noindent Proof: Let $\q$ be a $\s^d$-prime ideal of $R$. All we have to show that is that the fibre ring
\[(R\otimes_K K')\otimes_R k(\q)=k(\q)\otimes_K K'\] of $\q$ has non-empty $\s^\bullet$-spectrum. But this follows immediately from Theorem \ref{theo compatibility}. \qed

\subsection{Difference Kernels}

For the proof of the difference Chevalley theorem we need to generalize the theory of difference kernels (\cite[Chapter 6]{Cohn:difference}) from $\s$-fields to $\s$-pseudo fields.
This more or less straight forward task is realized in this section.

Let $K=e_1K\oplus\cdots\oplus e_dK$ be a $\s$-pseudo field. Note that if $R$ is a $K$-algebra then $R$ decomposes as a direct sum of rings $R=e_1R\oplus\cdots\oplus e_dR$. For every prime ideal $\q$ of $R$ there exists a unique $i\in\{1,\ldots,d\}$ such that $e_i\q$ is a prime ideal of $e_iR$ and $\q=e_1R\oplus\cdots\oplus e_i\q\oplus\cdots\oplus e_dR$.
If $R$ is a $K$-$\s$-algebra and $\q$ a $\s^{d'}$-prime ideal of $R$ then $d|d'$. Moreover a reflexive $\s$-ideal $\p$ of $R$ is a $\s$-pseudo prime ideal of $R$ with period $d$ if and only if $e_i\p$ is a prime ideal of $e_iR$ for $i=1,\ldots,d$.


By a \emph{difference kernel} over $K$ of length $t\geq 1$ in the $\s$-polynomial ring $K\{x\}=K\{x_1,\ldots,x_n\}$ we mean an ideal $\p_t$ of $K[x,\s(x),\ldots,\s^t(x)]$ such that
\begin{enumerate}
 \item the inverse image of $\p_t$ under $\s:K[x,\ldots,\s^{t-1}(x)]\rightarrow K[x,\ldots,\s^t(x)]$ equals $\p_{t-1}:=\p_t\cap K[x,\ldots,\s^{t-1}(x)]$ and
\item $e_i\p_t$ is a prime ideal of $e_iK[x,\ldots,\s^t(x)]$ for $i=1,\ldots,d$.
\end{enumerate}

If $\p$ is a $\s$-pseudo prime ideal of $K\{x\}$ of period $d$ and $t\geq 1$ then, obviously $\p_t=\p\cap K[x,\ldots,\s^t(x)]$ is a difference kernel of length $t$. Conversely, for a given difference kernel $\p_t$, a $\s$-pseudo prime ideal of $K\{x\}$ of period $d$ with $\p_t=\p\cap K[x,\ldots,\s^t(x)]$ is called a \emph{realization of $\p_t$}. By a \emph{prolongation} of $\p_t$ we mean a difference kernel $\p_{t+1}$ of length $t+1$ such that $\p_{t+1}\cap K[x,\ldots,\s^t(x)]=\p_t$.

From (i) we obtain an injection
\[ \s:K[x,\ldots,\s^{t-1}(x)]/\p_{t-1}\longrightarrow K[x,\ldots,\s^t(x)]/\p_t\]
that splits into injections \[\s:e_iK[x,\ldots,\s^{t-1}(x)]/e_i\p_{t-1}\longrightarrow e_{i+1}K[x,\ldots,\s^t(x)]/e_{i+1}\p_t.\]
Let $a,\ldots,\s^t(a)$ denote the image of $x,\ldots,\s^t(x)$ in $K[x,\ldots,\s^t(x)]/\p_t$. Then the above injections extend to injections of the quotient fields
\[\s:e_iK(e_ia,\ldots,e_i\s^{t-1}(a))\longrightarrow e_{i+1}K(e_{i+1}a,\ldots,e_{i+1}\s^{t-1}(a)).\]

Therefore we see that giving a difference kernel in $K\{x_1,\ldots,x_n\}$ is equivalent to specifying the following data: For $i=1,\ldots,d$ a field extension
$L_i=e_iK(a_i,\s(a_i),\ldots,\s^t(a_i))$ of $e_iK$ generated by certain $n$-tuples $a_i,\s(a_i),\ldots,\s^t(a_i)$ and field morphisms
$\s:L_i\rightarrow L_{i+1}$ such that $\s$ extends $\s:e_iK\rightarrow e_{i+1}K$ and $\s(\s^j(a_i))=\s^{j+1}(a_{i+1})$ for $j=0,\ldots,t-1$.

Set $L=L_1\oplus\cdots\oplus L_d$. Then we can recover $\p_t$ as the kernel of the $K$-morphism $K[x,\ldots,\s^t(x)]\rightarrow L$ defined by $e_i\s^j(x)\mapsto \s^j(a_i)$.

\begin{lemma}
 Every difference kernel $\p_t$ of length $t$ in $K\{x\}$ can be prolonged to a difference kernel $\p_{t+1}$ of length $t+1$.
\end{lemma}
\noindent Proof: We use the notation introduced above. For easier reading we abbreviate
\begin{align*}
 K_{i,t-1} & =e_iK\left(a_i,\ldots,\s^{t-1}(a_i)\right) \\
 K_{i+1,t} & =e_{i+1}K\left(a_{i+1},\ldots,\s^t(a_{i+1})\right) \\
 \widetilde{K}_{i+1} & =\s(e_iK)\left(\s(a_{i+1}),\ldots,\s^t(a_{i+1})\right)\subset K_{i+1,t} \\
\end{align*}
For $i=1,\ldots,d$ let $K_{i,t-1}[y]$ denote the polynomial ring in the variables $y=(y_1,\ldots,y_n)$ over the field $K_{i,t-1}$ and $\ida_i\subset K_{i,t-1}[y]$ the vanishing ideal of $\s^t(a_i)$. We have an isomorphism of fields
$\s: K_{i,t-1}\rightarrow \widetilde{K}_{i+1}$ that trivially extends to an isomorphism of polynomial rings $\psi: K_{i,t-1}[y]\rightarrow \widetilde{K}_{i+1}[y].$
We consider the base extension of $\psi(\ida_i)$ over $\widetilde{K}_{i+1}\subset K_{i+1,t}$. Let $\ida_i'\subset K_{i,t+1}[y]$ denote a minimal prime ideal above the ideal generated by $\psi(\ida_i)$ in $K_{i,t+1}[y]$. Then $\ida_i'\cap \widetilde{K}_{i+1}[y]=\psi(\ida_i)$. Let $\s^{t+1}(a_{i+1})$ denote the image of $y$ in the quotient field of $K_{i,t+1}[y]/\ida_i'$.
Passing to the quotient fields in
\[K_{i,t-1}[y]/\ida_i\rightarrow\widetilde{K}_{i+1}[y]/\psi(\ida_i)\hookrightarrow K_{i+1,t}[y]/\ida_i'\]
yields the desired extension of $\s$:
\[e_iK\left(a_i,\ldots,\s^t(a_i)\right)\rightarrow \sigma(e_iK)\left(\sigma(a_{i+1}),\ldots,\sigma^{t+1}(a_{i+1})\right)\rightarrow e_{i+1}K\left(a_{i+1},\ldots,\s^{t+1}(a_{i+1})\right).\]
\qed

\begin{cor} \label{cor realization of kernel}
Every difference kernel has a realization.
\end{cor}
\qed

\subsection{Inversive $\s$-rings}

In proofs it is often convenient to assume that certain $\s$-rings are inversive. This section provides the tools to make this reduction.
We recall that a $\s$-ring $R$ is called \emph{inversive} if $\s:R\rightarrow R$ is an automorphism.

\begin{defi}
Let $R$ be a $\s$-ring. By an \emph{inversive closure of $R$} we understand a pair $(u,R^*)$ where $R^*$ is an inversive $\s$-ring and $u:R\rightarrow R^*$ is a morphism of $\s$-rings satisfying the following universal property. For every morphism $f$ from $R$ into an inversive $\s$-ring $S$ there exists a unique morphism $g:R^*\rightarrow S$
making
\[
\begin{xy}
\xymatrix{
R  \ar[rr]^u \ar[dr]_-{f} & & R^* \ar@{.>}[ld]^g \\
& S &
}
\end{xy}
\]
commutative.
\end{defi}

As in \cite[Theorem 3, A.V.5]{Bourbaki:Algebra2} (see also \cite[Chapter 2, Section 5, p. 66]{Cohn:difference}) one sees that every $\s$-ring $R$ has an inversive closure $(u,R^*)$, which is then unique up to unique isomorphisms. Moreover
\[\ker(u)=\{r\in R|\ \exists\ n\geq 1 \text{ such that } \s^n(r)=0\}\] and for every $r\in R^*$ there exists $n\geq 0$ such that $\s^n(r)\in u(R)$.

Note that $u$ is injective if and only if $\s$ is injective. In this situation we will usually identify $R$ with the subring $u(R)$ of $R^*$ (and not mention $u$).
\begin{lemma} \label{lemma inversive homeom}
Let $R$ be a $\s$-ring and $(u,R^*)$ an inversive closure of $R$. Then $u:R\rightarrow R^*$ induces a bijection
$\s^\bullet\Spec(R^*)\rightarrow\s^\bullet\Spec(R),\ \q\mapsto u^{-1}(\q)$ that preserves inclusions and period.
For a $\s^d$-prime ideal $\q$ of $R$ the inverse is given by $\q^*=\{r^*\in R^*| \ \exists \ n\geq 0 \text{ such that } \s^{nd}(r^*)\in u(\q)\}$.
\end{lemma}
\noindent The simple proof is left to the reader.

%

\subsection{$\s$-algebraic Chevalley theorem}

The aim of this section is to prove Theorem \ref{theo difference algebraic Chevalley Introduction}. The overall plan of proof follows along the lines of the classical results of Cohn:
The problem breaks into two parts. First the case that $K=\Q(R)$ is relatively separably algebraically closed in $L=\Q(S)$ and secondly the case that $L$ is a separable algebraic field extension of $K$. By virtue of Babbitt's decomposition this second case can then be split into the case of a finite extension called the core of $L|K$ and the case of benign extensions.

\vspace{5mm}

We shall need a preparatory lemma. Let $R\subset S$ be an inclusion of integral domains and $a\subset S$ a finite generating set of $S$ as an $R$-algebra. Let $\q$ be a prime ideal of $R$.
Following \cite[Section 11, Chapter 7]{Cohn:difference} we call a prime ideal $\q'$ of $S$ with $\q'\cap R=\q$ a \emph{non-degenerate} lift of $\q$ to $S$ (with respect to $a$) if the following condition is satisfied: A subset of $a$ is algebraically independent over $\Q(R)$ (inside $\Q(S)$) if and only if the corresponding subset of the image of $a$ in the residue field $k(\q')$ is algebraically independent over $k(\q)$.

\begin{lemma} \label{lemma nondegenerate lift}
Let $R\subset S$ be an inclusion of integral domains and $a\subset S$ a finite generating set of $S$ as $R$-algebra. Assume that $\Q(R)$ is relatively separably algebraically closed in $\Q(S)$. Then there exists a non-zero $r\in R$ such that every prime ideal $\q$ of $R$ with $r\notin\q$ has a unique non-degenerate lift (with respect to $a$) to $S$.
\end{lemma}
\noindent Proof: This is \cite[Lemma 5, p. 209]{Cohn:difference}, only that there it is assumed that $R$ contains a field. But this is not necessary for the proof. The statement of  \cite[Lemma 3, p. 206]{Cohn:difference} (which is used in the proof of Lemma 5) without the assumption that $R$ contains a field can be found in \cite[Lemma 9.7.5, p. 77]{Grothendieck:EGAIV3}  \qed

\vspace{5mm}

The following proposition resolves the case when $K$ is relatively separably algebraically closed in $L$. The proposition and its proof are generalization of \cite[Theorem 4 (1), p. 211]{Cohn:difference}. The generalization is twofold: First we drop the assumption that $R$ is of finite $\s$-type over a $\s$-field (cf. \cite[Lemma 6.3.13, p. 390]{Levin:difference}). Secondly we provide uniformity in powers of $\sigma$.

\begin{prop}\label{prop Chev for primary}
Let $R\subset S$ be an inclusion of $\s$-domains such that $S$ is finitely $\s$-generated over $R$. Assume that $K=\Q(R)$ is relatively separably algebraically closed in $L=\Q(S)$. Then there exists a non-zero element $r\in R$ such that for all $d\geq 1$ every $\s^d$-prime ideal $\q$ of $R$ with $r\notin_\s\q$ has a lift to a $\sigma^d$-prime ideal $\q'$ of $S$.
\end{prop}
\noindent Proof: As the first step we show that we can assume without loss of generality that $R$ is inversive. Let $R^*\subset S^*$ denote the corresponding inclusion of
inversive closures. Then $R^*[S]\subset S^*$ is finitely $\s$-generated over $R^*$. Because $\Q(R)$ is relatively separably algebraically closed in $\Q(S)$ it follows that $\Q(R)^*$ is relatively separably algebraically closed in $\Q(S)^*$. Therefore $\Q(R^*)=\Q(R)^*$ is relatively separably algebraically closed in $\Q(R^*[S])\subset\Q(S^*)=\Q(S)^*$. In other words, the inclusion
$R^*\subset R^*[S]$ meets the requirements of the proposition. By assumption there exists a non-zero element $r^*\in R^*$ such that every $\s^d$-prime ideal $\q^*$ of $R^*$ with $r^*\notin_\s \q^*$ has a lift to a $\s^d$-prime ideal ${\q^*}'$ of $R^*[S]$. Let $n\geq 1$ be such that $r=\s^{n}(r^*)\in R$. Let $\q$ be a $\s^d$-prime ideal of $R$ with $r\notin_\s\q$.
As seen in Lemma \ref{lemma inversive homeom} the set $\q^*=\{s\in R^* | \ \exists \ m\geq 0 \text{ such that } \s^{md}(s)\in\q \}$ is a $\s^d$-prime ideal of $R^*$ with $\q^*\cap R=\q$. Suppose $r^*\in_\s\q^*$. Then there exists an integer $i\geq 0$ such that $\s^i(r^*)\in\q^*$, but then $\s^{md+i}(r^*)\in\q$ for all $m$ greater than or equal to some $m'\geq 0$. This contradicts $r=\s^n(r^*)\notin_\s\q$. Thus $r^*\notin_\s\q^*$ and so there exists a $\s^d$-prime ideal ${\q^*}'$ of $R^*[S]$ with ${\q^*}'\cap R^*=\q^*$. Now $\q'={\q^*}'\cap S$ has the desired properties.

\vspace{5mm}

From now on we assume that $R$ is inversive.
Because $S$ is finitely $\s$-generated over $R$ we can write $S=R\{a\}$ for some $a=(a_1,\ldots,a_n)\in S^n$. Let $\ida'\subset K\{x\}=K\{x_1,\ldots,x_n\}$ denote the $\s$-prime ideal of all $\s$-polynomials vanishing on $a$. Set $\ida=R\{x\}\cap \ida'$. Then $\ida$ is a $\s$-prime ideal of $R\{x\}$ and $S=R\{x\}/\ida$. Let $\{p_1,\ldots,p_m\}\subset\ida'$ be a characteristic set of $\ida'$ (see \cite[Section 2.4]{Levin:difference}). (We use the standard ranking, i.e. $x_1<\cdots<x_n<\s(x_1)<\cdots$ .) If we multiply $p_1,\ldots,p_m$ with a non-zero $\lambda\in K$ then $\lambda p_1,\ldots,\lambda p_m$ is again a characteristic set of $\ida'$. Thus we can assume that $p_1,\ldots,p_m\in\ida\subset R\{x\}$.
Let $q\in R\{x\}$ denote the product of the initials of $p_1,\ldots,p_m$ and let $t\geq 1$ be an integer such that the order of each $p_i$ is at most $t$.

Consider the inclusion $R\subset R[a,\s(a),\ldots,\s^{t}(a)]$. Because $\Q(R)$ is relatively separably algebraically closed in $\Q(S)$ we also know that $\Q(R)$ is relatively separably algebraically closed in $\Q(R[a,\s(a),\ldots,\s^{t}(a)])$. Note that $p_1(a)=0,\ldots,p_m(a)=0$.

By \cite[Proposition 2.4.4, p. 132]{Levin:difference} we know that $q\notin\ida'$. Thus $q\notin\ida$ and consequently $q(a)\in R[a,\s(a),\ldots,\s^{t}(a)]$ is non-zero. It follows from Lemma \ref{lemma nondegenerate lift} that there exists a non-zero $r\in R$ such that for every prime ideal $\q$ of $R$ with $r\notin\q$ there exists a prime ideal $\widehat{\q}$ of $R[a,\s(a),\ldots,\s^{t}(a)]$ with the following properties:
\begin{enumerate}
 \item $\widehat{\q}$ is a non-degenerate lift (with respect to $a,\s(a),\ldots,\s^{t}(a)$) of $\q$ to $R[a,\ldots,\s^{t}(a)]$,
\item $\widehat{\q}\cap R[a,\ldots,\s^{t-1}(a)]$ is the unique non-degenerate lift of $\q$ to $R[a,\ldots,\s^{t-1}(a)]$ and
\item $q(a)\notin \widehat{\q}$.
\end{enumerate}

Now let $\q$ be a $\s^d$-prime ideal of $R$ with $r\notin_\s\q$. Set $\q_1=\q=\s^{-d}(\q),\q_2=\s^{-(d-1)}(\q),\ldots,\q_d=\s^{-1}(\q)$. Then $r\notin\q_1,\ldots,\q_d$ and $\s:R\rightarrow R$ induces $\s:k(\q_i)\rightarrow k(\q_{i+1})$.
Let $i\in\{1,\ldots,d\}$ and let $\widetilde{\q_i}$ denote the prime ideal of $R[a,\ldots,\s^{t}(a)]\otimes_R k(\q_i)$ corresponding to $\widehat{\q_i}$ (where $\widehat{\q_i}$ is as specified above). Because $R$ is inversive we have an isomorphism
\begin{align*}
\psi_i: R[a,\ldots,\s^{t-1}(a)]\otimes_R k(\q_i) & \longrightarrow R[\s(a),\ldots,\s^t(a)]\otimes_R k(\q_{i+1}) \\
 f\otimes g & \longmapsto \s(f)\otimes\s(g).
\end{align*}
Let
\[\alpha_i: R[a,\ldots,\s^{t-1}(a)]\otimes_R k(\q_i)\longrightarrow R[a,\ldots,\s^{t}(a)]\otimes_R k(\q_i)\]
and
\[\beta_{i+1}: R[\s(a),\ldots,\s^{t}(a)]\otimes_R k(\q_{i+1})\longrightarrow R[a,\ldots,\s^{t}(a)]\otimes_R k(\q_{i+1})\]
denote the morphisms corresponding to the inclusions of rings. Then $\beta_{i+1}^{-1}(\widetilde{\q_{i+1}})$ corresponds to a non-degenerate lift of $\q_{i+1}$ to
$R[\s(a),\ldots,\s^t(a)]$. Because $\psi_i$ is an isomorphism, $\psi_i^{-1}(\beta_{i+1}^{-1}(\widetilde{\q_{i+1}}))$ corresponds to a non-degenerate lift of $\q_i$ to $R[a,\ldots,\s^{t-1}(a)]$. But also $\alpha_i^{-1}(\widetilde{\q_i})$ provides a non-degenerate lift of $\q$ to $R[a,\ldots,\s^{t-1}(a)]$. By the uniqueness of such a lift we conclude that $\alpha_i^{-1}(\widetilde{\q_i})=\psi_i^{-1}(\beta_{i+1}^{-1}(\widetilde{\q_{i+1}}))$.

Let $a_i,\ldots,\s^t(a_i)$ denote the images of $a,\ldots,\s^t(a)$ in $k(\widetilde{\q_i})=k(\widehat{\q_i})$. Then
$\alpha_i^{-1}(\widetilde{\q_i})=\psi_i^{-1}(\beta_{i+1}^{-1}(\widetilde{\q_{i+1}}))$ signifies that we have well defined morphisms
\[\s:k(\q_i)\left(a_i,\ldots,\s^{t-1}(a_i)\right)\longrightarrow k(\q_{i+1})\left(a_{i+1},\ldots,\s^{t}(a_{i+1})\right)\]
extending $\s:k(\q_i)\rightarrow k(\q_{i+1})$ given by $\s(a_i)=\s(a_{i+1}),\ldots,\s(\s^{t-1}(a_i))=\s^t(a_{i+1})$.
We have thus constructed a difference kernel $\widetilde{\p}_t\subset K[x,\ldots,\s^t(x)]$ of length $t$ over the $\s$-pseudo field $K=k(\q_1)\oplus\cdots\oplus k(\q_d)$.
By Corollary \ref{cor realization of kernel} there exists a $\s$-pseudo prime ideal $\widetilde{\p}$ of $K\{x\}$ with period $d$ such that $\widetilde{\p}\cap K[x,\ldots,\s^t(x)]=\widetilde{\p}_t$.

Let $\p=\q_1\cap\cdots\cap\q_d$ denote the $\s$-pseudo prime ideal of $R$ associated with $\q$. The residue map $R\rightarrow  k(\p)=K$ extends to a morphism
$\phi: R\{x\}\rightarrow k(\widetilde{\p})$ by $\phi(x)=a_1+\cdots+a_d$. Let $\p''$ denote the kernel of $\phi$. Then $\p''$ is a $\s$-pseudo prime ideal of $R\{x\}$ of period $d$ with $\p''\cap R=\p$.

We will show that $\ida\subset \p''$. By construction $p_1,\ldots,p_m\in\p''$ and $q\notin_\s\p''$.
Let $p\in\ida$. By \cite[Proposition 2.4.4, p. 132]{Levin:difference} there exist integers $n_1,\ldots,n_r\geq 0$ and $m_1,\ldots,m_r\geq 1$ such that $\s^{n_1}(q)^{m_1}\cdots\s^{n_r}(q)^{m_r}p$ lies in the $\s$-ideal generated by $p_1,\ldots,p_m$, in particular $\s^{n_1}(q)^{m_1}\cdots\s^{n_r}(q)^{m_r}p\in\p''$. As $q\notin_\s\p''$ this implies $p\in\p''$. Therefore $\p''$ contains $\ida$ and defines a $\s$-pseudo prime ideal $\p'$ of $S=R\{x\}/\ida$ with period $d$ such that $\p'\cap R=\p$. One of the minimal prime ideals of $\p'$, call it $\q'$, satisfies $\q'\cap R=\q$. \qed

\begin{lemma} \label{lemma chevalley core}
Let $R\subset S$ be an inclusion of $\s$-domains such that $S$ is finitely generated over $R$ (not merely finitely $\sigma$-generated). Assume that $L=\Q(S)$ is an algebraic extension of $K=\Q(R)$ and let $s$ be a non-zero element of $S$. Then there exists an integer $l\geq 1$ and a non-zero element $r$ of $R$ such that for all $d\geq 1$ every $\s^d$-prime ideal $\q$ of $R$ with $r\notin_\s\q$ has a lift to a $\s^{ld}$-prime ideal $\q'$ of $S$ with $s\notin_\s\q'$.
\end{lemma}
\noindent Proof:
As the first step we will show that we can assume without loss of generality that $R$ is inversive. Let $R^*\subset S^*$ be the inclusion of the inversive closures. Then
$R^*[S]\subset S^*$ is finitely generated as $R^*$-algebra. Because $\Q(S)$ is algebraic over $\Q(R)$ it follows that $\Q(S)^*$ is algebraic over $\Q(R)^*$. Therefore
$\Q(R^*[S])\subset\Q(S^*)=\Q(S)^*$ is algebraic over $\Q(R)^*=\Q(R^*)$. So the inclusion $R^*\subset R^*[S]$ satisfies the assumptions of the lemma. Thus, by assumption
there exists an integer $l\geq 1$ and a non-zero $r^*\in R^*$ such that every $\s^d$-prime ideal $\q^*$ of $R^*$ with $r^*\notin_\s\q^*$ lifts to a $\s^{ld}$-prime ideal ${\q^*}'$ of $R^*[S]$ with $s^*=s\notin_\s{\q^*}'$. There is an integer $n\geq 0$ such that $r=\s^n(r^*)$ lies in $R$. Let $\q$ be a $\s^d$-prime
ideal of $R$ with $r\notin_\s\q$. From Lemma \ref{lemma inversive homeom} it follows that $\q$ is of the form $\q=\q^*\cap R$ with
$\q^*=\{s\in R^*| \ \exists \ m\geq 0 \text{ such that } \s^{md}(s)\in \q\}$. As $r\notin_\s \q$ implies $r^*\notin_\s\q^*$ there exists a $\s^{ld}$-prime ideal ${\q^*}'\subset R^*[S]$ with $s\notin_\s {\q^*}'$ and ${\q^*}'\cap R^*=\q^*$. Consequently $\q'={\q^*}'\cap S$ is a $\s^{ld}$-prime ideal of $S$ with $\q'\cap R=\q$ and $s\notin_\s\q'$ as desired.

\vspace{5mm}

From now on we will assume that $R$ is inversive. It follows from the assumptions that $L=K[S]$ is a finite field extension of $K$. Because $K$ is inversive this implies that $L=K[S]$ is also inversive. Let $a=(a_1,\ldots,a_n)\in S^n$ such that $S=R[a]$. Then $a$ generates $L$ as an $K$-algebra. Because $L$ is inversive also $\s(a)$ generates $L$ as $K$-algebra. Therefore every $a_i$ is of the form $a_i=P_i(\s(a))$ for some $P_i\in K[x]=K[x_1,\ldots,x_n]$. We can choose $r_1\in R$ such
that $P_i\in R_{r_1}[x]$ for $i=1,\ldots,n$.

The elements $a_1,\ldots,a_n,\frac{1}{s}\in L$ are all algebraic over $K$. Thus there exists a non-zero $r_2\in R$ such that $a_1,\ldots,a_n,\frac{1}{s}$ are integral
over $R_{r_2}$.

Set $r=r_1r_2\in R$ and let $T$ denote the multiplicatively closed subset of $R$ generated by $\s^i(r)$ for $i\in\mathbb{Z}$. Then $T^{-1}R\subset L$ is an inversive $\s$-ring, and $T^{-1}R[S]\subset L$ is also inversive because by construction every $a_i$ lies in the image of $\s:T^{-1}R[S]\rightarrow T^{-1}R[S]$.
Moreover $(T^{-1}R)[S]\{\frac{1}{s}\}$ is integral over $T^{-1}R$.
As $T^{-1}R[S]$ is integral and finitely generated over $T^{-1}R$ it follows that $T^{-1}R[S]$ is finitely generated as $T^{-1}R$-module. Let $l'$ denote the cardinality of a finite generating set of $T^{-1}R[S]$ as $T^{-1}R$-module and set $l=l'!$.

Let $\q$ be a $\s^d$-prime ideal of $R$ with $r\notin_\s\q$. Then $\q$ lifts to a $\s^d$-prime ideal $\widetilde{\q}$ of $T^{-1}R$.
Because $(T^{-1}R)[S]\{\frac{1}{s}\}$ is integral over $T^{-1}R$ every prime ideal of $T^{-1}R$ lifts to a prime ideal of $(T^{-1}R)[S]\{\frac{1}{s}\}$ (see e.g. \cite[Theorem 5.10, p. 62]{AtiyahMacdonald:Introductiontocommutativealgebra}).
In particular $\widetilde{\q}$ lifts to a prime ideal $\widetilde{\q}'$ of $T^{-1}R[S]$ with $s\notin_\s\widetilde{\q}'$.
Now it only remains to see that $\s^{-ld}(\widetilde{\q}')=\widetilde{\q}'$ because then
$\q'=S\cap\widetilde{\q}'$ is the desired lift of $\q$ to $S$.

The fibre ring \[\operatorname{FR}(\widetilde{\q})=T^{-1}R[S]\otimes_{T^{-1}R} k(\widetilde{\q})\] of $T^{-1}R\subset T^{-1}R[S]$ over $\widetilde{\q}$ is not the zero ring and generated as $k(\widetilde{\q})$-vector space by $l'$ elements. This implies that $\operatorname{FR}(\widetilde{\q})$ has at most $l'$ prime ideals. Because all the rings involved in the definition of
$\operatorname{FR}(\widetilde{\q})$ are inversive also $\operatorname{FR}(\widetilde{\q})$ is inversive. Thus $\overline{\q}\mapsto\s^{-d}(\overline{\q})$ defines a permutation of the prime ideals of $\operatorname{FR}(\widetilde{\q})$. In particular, the prime ideal $\overline{\q}$ of $\operatorname{FR}(\widetilde{\q})$ corresponding to $\widetilde{\q}'\subset T^{-1}R[S]$ satisfies
$\s^{-\overline{l}d}(\overline{\q})=\overline{\q}$ for some $\overline{l}\leq l'$. But then also $\s^{-ld}(\overline{\q})=\overline{\q}$ and $\s^{-ld}(\widetilde{\q}')=\widetilde{\q}'$ as desired. \qed

\vspace{5mm}

The following lemma and its proof are a generalization of part of \cite[Theorem 6, p. 217]{Cohn:difference}.

\begin{lemma} \label{lemma benign}
Let $R$ be a $\s$-domain and set $K=\Q(R)$. Let $L=K\langle a\rangle$ be a benign $\s$-field extension of $K$ with
minimal standard generator $a\in L$ and let $s$ be a non-zero element of $S=R\{a\}$. Then there exists a non-zero $r\in R$ such that every $\s^d$-prime ideal $\q$ of $R$ with
$r\notin_\s\q$ lifts to a $\s^d$-prime ideal $\q'$ of $S$ with $s\notin_\s\q'$.
\end{lemma}
\noindent Proof: Let $f',g'\in K[x]$ denote the minimal polynomials of $a$ and $s$ respectively over $K$. Then there exists a non-zero $r_1\in R$ such that $f=r_1f'$ and $g=r_1g'$ have coefficients in $R$. Let $r_2=g(0)$ denote the trailing coefficient of $g$ and set $r=r_1r_2$.

Let $\q$ be a $\s^d$-prime ideal of $R$ with $r\notin_\s \q$ and let $\p=\q\cap\cdots\cap\s^{-(d-1)}(\q)$ denote the $\s$-pseudo prime ideal of $R$ associated with $\q$.
The residue map $R\rightarrow k(\p)$ trivially extends to $R\{x\}\rightarrow k(\p)\{x\},\ p\mapsto \overline{p}$.

We claim that $\overline{f}$ has a solution $b\in M$ in a $\s$-pseudo field extension $M$ of $k(\p)$ of period $d$. This is easy to see: If $k(\p)=e_1k(\p)\oplus\cdots\oplus e_dk(\p)$, let $M_1,\ldots,M_d$ denote algebraic closures of $e_1k(\p),\ldots,e_dk(\p)$. Then the injections $\s:e_ik(\p)\rightarrow e_{i+1}k(\p)$ extend to $\s:M_i\rightarrow M_{i+1}$, thereby defining a $\s$-pseudo field extension $M=M_1\oplus\cdots\oplus M_d$ of $k(\p)$. By construction $e_i\overline{f}\in M_i[x]$ is a polynomial of positive degree for $i=1,\ldots,d$. We conclude that $\overline{f}$ has a solution $b$ in $M$.

Let $\phi:R\{x\}\rightarrow M$ denote the extension of $R\rightarrow k(\p)$ defined by $x\mapsto b$. Let $\widetilde{\p}$ denote the kernel of $\phi$. Then $\widetilde{\p}$ is a $\s$-pseudo prime ideal of $R\{x\}$ of period $d$ with $\widetilde{\p}\cap R=\p$.

Let $\ida'\subset K\{x\}$ denote the $\s$-prime ideal of $\s$-polynomials vanishing on $a$. Because $L$ is benign over $K$ with minimal standard generator $a$ one immediately sees that $f'$ is a characteristic set of $\ida'$. Set $\ida=R\{x\}\cap\ida'$. We will show that $\ida\subset\widetilde{\p}$. Obviously the $\s$-ideal generated by $f$ in $R\{x\}$ lies in $\widetilde{\p}$. Let $p\in\ida\subset\ida'$. Note that $f=r_1f'$ also is a characteristic set of $\ida'$ and that the initial of $f$ is $r_1$. If follows from \cite[Proposition 2.4.4, p. 132]{Levin:difference} that there exist integers $n_1,\ldots,n_k\geq 0$, $m_1,\ldots,m_k\geq 1$ such that $\s^{n_1}(r_1)^{m_1}\cdots\s^{n_k}(r_1)^{m_k}p$ lies in the $\s$-ideal generated by $f$. In particular
$\s^{n_1}(r_1)^{m_1}\cdots\s^{n_k}(r_1)^{m_k}p\in\widetilde{\p}$. By construction $r_1\notin_\s\widetilde{\p}$, therefore $p\in\widetilde{\p}$. Thus $\ida\subset\widetilde{\p}$ and $\widetilde{\p}$ defines a $\s$-pseudo prime ideal $\p'$ of period $d$ in $S=R\{a\}=R\{x\}/\ida$ with $\p'\cap R=\p$. We have $\q'\cap R=\q$ for a minimal prime ideal $\q'$ of $\p'$.

It remains to see that $s\notin_\s\q$. Suppose $\s^i(s)\in \q'$ for some $i\geq 0$. Because $g(s)=0$ this implies $\s^i(r_2)\in\q'\cap R=\q$. This contradicts $r=r_1r_2\notin_\s\q$. \qed

\vspace{5mm}

\noindent The following proposition resolves the separable algebraic part.

\begin{prop} \label{prop chevalley separabel algebraic part}
Let $R\subset S$ be an inclusion of $\s$-domains such that $S$ is finitely $\s$-generated over $R$. Assume that $L=\Q(S)$ is a separable algebraic extension of $K=\Q(R)$ and let $s$ be a non-zero element of $S$. Then there exists an integer $l\geq 1$ and a non-zero element $r\in R$ such that every $\s^d$-prime ideal $\q$ of $R$ with $r\notin_\s\q$ has a lift to a $\s^{ld}$-prime ideal $\q'$ of $S$ with $s\notin_\s\q'$.
\end{prop}
\noindent Proof: As in Proposition \ref{prop Chev for primary} or Lemma \ref{lemma chevalley core} one shows that we can assume without loss of generality that $R$ is inversive. Moreover, replacing $S$ with $S\{\frac{1}{s}\}$ we can assume $s=1$, i.e., we can neglect $s$.
%
%
The $\s$-field extension $L|K$ is finitely $\s$-generated and separable algebraic.
Let $\widetilde{L}$ denote a normal closure of $L|K$. By this we mean a $\s$-field extension of $L$ such that $\widetilde{L}$ is the field theoretic normal closure of $L|K$ .
Such an extension exists and is finitely $\s$-generated (see \cite[Section 16, p. 216]{Cohn:difference}). Because $L|K$ is separable, $\widetilde{L}|K$ is also separable.
Thus $\widetilde{L}|K$ is a normal separable algebraic finitely $\s$-generated $\s$-field extension of the inversive $\s$-field $K$, i.e., $\widetilde{L}|K$ meets the requirements for Babbitt's decomposition (\cite[Theorem 5.4.13, p. 336 ]{Levin:difference}). Thus there exist $a_1,\ldots,a_n\in\widetilde{L}$ such that $\widetilde{L}=\widetilde{L}_K\langle a_1,\ldots,a_n\rangle$ and $\widetilde{L}_K\langle a_1,\ldots,a_i\rangle$ is a benign $\s$-field extension of $\widetilde{L}_K\langle a_1,\ldots,a_{i-1}\rangle$ with minimal standard generator $a_i$ for $i=1,\ldots,n$.

The core $\widetilde{L}_K$ of $\widetilde{L}$ over $K$ is a finite separable algebraic extension of $K$. Thus there exists a primitive element $a_0\in\widetilde{L}_K$, i.e., $\widetilde{L}_K=K[a_0]$. Let $f\in K[x]$ such that $\s(a)=f(a)$. Let $r_1\in R$ such that $f$ has coefficients in $R_{r_1}$ and set $\widetilde{R}=R\{\frac{1}{r_1}\}$.
Then $\widetilde{R}[a_0]\subset\widetilde{L}$ is a difference ring.

As $\Q(S)=L\subset \widetilde{L}=\Q(\widetilde{R}\{a_0,\ldots,a_n\})$ and $S$ is finitely $\s$-generated over $R$ there exists a non-zero $r_n\in\widetilde{R}\{a_0,\ldots,a_n\}$ such that $S\subset \widetilde{R}\{a_0,\ldots,a_n,\frac{1}{r_n}\}$.

Consider the inclusion $\widetilde{R}\{a_0,\ldots,a_{n-1}\}\subset\widetilde{R}\{a_0,\ldots,a_n\}$ of $\s$-domains. It follows from Lemma \ref{lemma benign} that there exists a non-zero $r_{n-1}\in\widetilde{R}\{a_0,\ldots,a_{n-1}\}$ such that every $\s^d$-prime ideal $\q_{n-1}$ of $\widetilde{R}\{a_0,\ldots,a_{n-1}\}$ with $r_{n-1}\notin_\s\q_{n-1}$ lifts to a $\s^d$-prime ideal $\q_n$ of $\widetilde{R}\{a_0,\ldots,a_n\}$ with $r_n\notin_\s\q_n$.
Applying Lemma \ref{lemma benign} to the inclusion $\widetilde{R}\{a_0,\ldots,a_{n-2}\}\subset\widetilde{R}\{a_0,\ldots,a_{n-1}\}$ we see that there exists a non-zero $r_{n-2}\in \widetilde{R}\{a_0,\ldots,a_{n-2}\}$ such that
every $\s^d$-prime ideal $\q_{n-2}$ of $\widetilde{R}\{a_0,\ldots,a_{n-2}\}$ with $r_{n-2}\notin_\s\q_{n-2}$ lifts to a $\s^d$-prime ideal of $\widetilde{R}\{a_0,\ldots,a_{n-1}\}$ with $r_{n-1}\notin_\s\q_{n-1}$. Combining the above two statements yields that every $\s^d$-prime ideal $\q_{n-2}$ of $\widetilde{R}\{a_0,\ldots,a_{n-2}\}$ with $r_{n-2}\notin_\s\q_{n-2}$ lifts to a $\s^d$-prime ideal $\q_n$ of $\widetilde{R}\{a_0,\ldots,a_n\}$ with $r_n\notin_\s\q_n$. Inductively we find that there exists a non-zero $r_0\in \widetilde{R}\{a_0\}=\widetilde{R}[a_0]$ such that every $\s^d$-prime ideal $\q_0$ of $\widetilde{R}\{a_0\}$ with $r_0\notin_\s\q_0$ lifts to a $\s^d$-prime ideal $\q_n$ of $\widetilde{R}\{a_0,\ldots,a_n\}$ with $r_n\notin_\s\q_n$.
Applying Lemma \ref{lemma chevalley core} to the inclusion $\widetilde{R}\subset\widetilde{R}[a_0]$ yields a non-zero element $\widetilde{r}\in\widetilde{R}$ and an integer $l\geq 1$ such that every $\s^d$-prime ideal $\widetilde{\q}$ of $\widetilde{R}$ with $\widetilde{r}\notin_\s\widetilde{\q}$ lifts to a $\s^{ld}$-prime ideal $\q_0$ of $\widetilde{R}\{a_0\}$ with $r_0\notin_\s\q_0$. Therefore every $\s^d$-prime ideal $\widetilde{\q}$ of $\widetilde{R}$ with $\widetilde{r}\notin_\s\widetilde{\q}$ lifts to a $\s^{ld}$-prime ideal $\q_n$ of $\widetilde{R}\{a_0,\ldots,a_n\}$ with $r_n\notin_\s\q_n$. As $\widetilde{R}=R\{\frac{1}{r_1}\}$ there is a non-zero $r\in R$ such that every $\s^d$-prime ideal $\q$ of $R$ with $r\notin_\s\q$ lifts to a $\s^d$-prime ideal $\widetilde{\q}$ of $\widetilde{R}$ with $\widetilde{r}\notin_\s\widetilde{\q}$. Therefore every $\s^d$-prime ideal $\q$ of $R$ with $r\notin_\s\q$ lifts to a $\s^{ld}$-prime ideal $\q_n$ of $\widetilde{R}\{a_0,\ldots,a_n\}$ with $r_n\notin_\s \q_n$. Thus $\q$ lifts to a $\s^{ld}$-prime ideal $\widehat{\q}$ of
$\widetilde{R}\{a_0,\ldots,a_n,\frac{1}{r_n}\}$. Finally $\q'=S\cap \widehat{\q}$ is the desired lift of $\q$.\qed

\vspace{5mm}

It is now a simple matter to combine Propositions \ref{prop Chev for primary} and \ref{prop chevalley separabel algebraic part} to obtain the general theorem.

\begin{theo}[$\s$-algebraic Chevalley theorem] \label{theo difference algebraic Chevalley}
Let $R\subset S$ be an inclusion of $\s$-domains such that $S$ is finitely $\s$-generated over $R$. Then there exists an integer
$l\geq 1$ and a non-zero element $r\in R$ such that for every $\s^d$-prime ideal $\q$ of $R$ with $r\notin_\s\q$ there exists a $\sigma^{ld}$-prime ideal $\q'$ of $S$ with $\q'\cap R=\q$.
\end{theo}
\noindent Proof: Set $K=\Q(R)$ and $L=\Q(S)$. Then $L$ is a finitely $\s$-generated $\s$-field extension of $K$. Let $\widetilde{K}\subset L$ denote the relative separable algebraic closure of $K$ in $L$.
Then $\widetilde{K}$ is a finitely $\s$-generated $\s$-extension of $K$ (see \cite[Theorem 18, p. 145]{Cohn:difference}), say $\widetilde{K}=K\langle a\rangle$ with $a\in \widetilde{K}^n$.
Applying Proposition \ref{prop Chev for primary} to the inclusion $R\{a\}\subset R\{a\}[S]$ yields a non-zero $\widetilde{r}\in R\{a\}$ such that every $\s^d$-prime ideal $\widetilde{\q}$ of $R\{a\}$ with $\widetilde{r}\notin_\s\widetilde{\q}$ lifts to a $\s^d$-prime ideal $\widehat{\q}$ of $R\{a\}[S]$.

Applying Proposition \ref{prop chevalley separabel algebraic part} to the inclusion $R\subset R\{a\}$ yields a non-zero $r\in R$ and an integer $l\geq 1$ such that every $\s^d$-prime ideal $\q$ of $R$ with $r\notin_\s\q$ lifts to a $\s^{ld}$-prime ideal $\widetilde{\q}$ of $R\{a\}$ with $\widetilde{r}\notin_\s\widetilde{\q}$. Combining the above two statements we find that every $\s^d$-prime ideal $\q$ of $R$ with $r\notin_\s\q$ lifts to a $\s^{ld}$-prime ideal $\widehat{\q}$ of $R\{a\}[S]$. Finally $\q'=\widehat{\q}\cap S$ is the desired lift of $\q$. \qed

\vspace{5mm}

We pause to explain that Theorem \ref{theo difference algebraic Chevalley} above is more or less optimal. It was already observed in the introduction (Example \ref{ex failure of naive chevalley}) that one can not in general expect to have $l=1$. Looking at an inclusion of the form $R\subset R\{\frac{1}{r}\}$ will convince the reader that the condition $r\notin_\s\q$ can not be relaxed to $r\notin\q$ (and is maybe more natural).

Let $R\subset S$ be an inclusion of $\s$-domains such that $S$ is finitely $\s$-generated over $R$. One might hope that there exists an integer $d\geq 1$ such that the ``naive'' version of a $\s$-algebraic Chevalley theorem holds for $(R,\s^d)\subset (S,\s^d)$, i.e., one might hope for the validity of the following variation:

{\bf Statement A:} There exist an integer $d\geq 1$ and a non-zero element $r\in R$ such that for every $\s^d$-prime ideal $\q$ of $R$ with $r\notin\q$ there exists a $\s^d$-prime ideal $\q'$ of $S$ with $\q=\q'\cap R$.

\vspace{3mm}

However Statement A is NOT true. The following example provides a counterexample.

\begin{ex}
Let $q\in\mathbb{C}^\times$ and choose a root $\sqrt{q}\in\mathbb{C}^\times$, i.e., $\sqrt{q}^2=q$. Let $R=\mathbb{C}[t,t_1,t_2,\ldots ]$ denote the polynomial ring in  the indeterminates $t,t_1,t_2,\ldots$. For $i\geq 1$ let $c_i\in\mathbb{C}^\times$ denote a complex number such that $c_i^i=-\sqrt{q}^i$. We endow $R$ with the structure of a $\s$-ring by letting $\s$ act on $\mathbb{C}$ as the identity and setting $\s(t)=qt$ and $\s(t_i)=c_it_i$ for $i\geq 1$.
Set $S=R[\sqrt{t}]$ with $\s(\sqrt{t})=\sqrt{q}\sqrt{t}$.

Suppose that Statement A is true. Then there exists an integer $d\geq 1$ such that every $\s^d$-prime ideal $\q$ of $R$ with $r\notin\q$ lifts to a $\s^d$-prime ideal $\q'$ of $S$.

Let $\lambda\in\mathbb{C}$. Then
\[\s^d(t-\lambda^2t_d^2)=q^dt-\lambda^2(c_dt_d)^2=q^dt-\lambda^2q^dt_d^2=q^d(t-\lambda^2t_d^2).\]

This shows that $\q_\lambda=(t-\lambda^2t_d^2)$ is a $\s^d$-prime ideal of $R$.
Suppose $\q_\lambda$ lifts to a $\s^d$-prime ideal $\q_\lambda'$ of $S$. For $s\in S$ we let $\overline{s}$ denote the image of $s\in $ in $k(\q_\lambda')$. Then
$\s^d(\lambda\overline{t_d})=\lambda c_d^d\overline{t_d}=-\sqrt{q}^d(\lambda \overline{t_d})$. On the other hand
$\s^d\left(\overline{\sqrt{t}}\right)=\sqrt{q}^d\overline{\sqrt{t}}$. Because both $\overline{\sqrt{t}}$ and $\lambda\overline{t_d}$ are roots of $\overline{t}$ this is impossible. I.e. $\q_\lambda$ does not lift to a $\s^d$-prime ideal of $S$.

But then $r$ must lie in every $\q_\lambda$ which is impossible. \qed
\end{ex}

\vspace{5mm}

We will conclude this section with some reformulations of Theorem \ref{theo difference algebraic Chevalley}.
Akin to Theorem \ref{theo usual algebraic chevalley with closed field} we can also give a formulation of Theorem \ref{theo difference algebraic Chevalley} with $\s$-closed fields, i.e., models of ACFA. (We recall that a $\s$-field $k$ is a model of ACFA if for every finitely $\s$-generated $k$-$\s$-domain $R$ there exists a $k$-$\s$-morphism $R\rightarrow k$.) However, because of possibly occurring incompatibility the statement is less uniform in powers of $\s$.

\begin{cor}[$\s$-algebraic Chevalley theorem with $\sigma$-closed fields]
Let $R\subset S$ be an inclusion of $\s$-domains such that $S$ is finitely $\s$-generated over $R$. Then there exists a non-zero element $r\in R$ such that every $\s^d$-morphism $\psi: R\rightarrow k$ into a $\s$-closed field $k$ with $\psi(r\s(r)\cdots\s^{d-1}(r))\neq 0$ can be extended to a $\s^{ld}$-morphism $\widetilde{\psi}:S\rightarrow k$ (where $l\geq 1$ is allowed to depend on $\psi$).
\end{cor}
\noindent Proof: By Theorem \ref{theo difference algebraic Chevalley} there exists an integer $l'\geq 1$ and a non-zero $r\in R$ such that every $\s^d$-prime ideal $\q$ of $R$ with $r\notin_\s\q$ lifts to a $\s^{l'd}$-prime ideal $\q'$ of $S$.
Let $\psi: R\rightarrow k$ be a $\s^d$-morphism into a $\s$-closed field $k$ with $\psi(r\s(r)\cdots\s^{d-1}(r))\neq 0$. Then the kernel $\q$ of $\psi$ is a $\s^d$-prime ideal of $R$ with $r\notin_\s\q$. Thus there exists a $\s^{l'd}$-prime ideal $\q'$ of $S$ with $\q'\cap R=\q$. The $\s^{l'd}$-field extension $k(\q')|k(\q)$ is finitely $\s^{l'd}$-generated. Thus by Corollary \ref{cor compatibility} there exists an integer $d'$ such that $k(\q')|k(\q)$ and $k|k(\q)$ are compatible as $\s^{d'l'd}$-fields. Set $l=d'l'$ and let $M$ be a $\s^{ld}$-field containing $k(\q')$ and $k$. Then $k\{S/\q'\}_{\s^{ld}}=k[S/\q']\subset M$ is a finitely $\s^{ld}$-generated $k$-$\s^{ld}$-domain. As $(k,\s^{ld})$ also is a model of ACFA (see \cite[Corollary 1.12, p. 3013]{Hrushovskietal:ModelTheoryofDifferencefields}), we see that there exists a $k$-$\s^{ld}$-morphism $k[S/\q']\rightarrow k$. Then
\[\widetilde{\psi}:S\rightarrow k[S/\q']\rightarrow k\] is the desired extension of $\psi$. \qed

\vspace{5mm}

The most useful form of Theorem \ref{theo difference algebraic Chevalley} for the applications in the second part is the following.

\begin{cor}\label{cor chevalley for pseudo domain}
Let $R\subset S$ be an inclusion of $\s$-pseudo domains such that $S$ is finitely $\s$-generated over $R$. Then there exists a non-zero divisor $r\in R$ such that every $\s$-pseudo prime ideal $\p$ of R with $r\notin_\s\p$ lifts to a $\s$-pseudo prime ideal of $S$.
\end{cor}
\noindent Proof: Let $\ida_1,\ldots,\ida_n$ denote the minimal prime ideals of $R$. Then for $i=1,\ldots,n$ there is a minimal prime ideal ${\ida_i}'$ of $S$ with ${\ida_i}'\cap R=\ida_i$. Let $n'$ denote the period of ${\ida_i}'$ (which does not depend on $i$ because $S$ is a $\s$-domain). For $i=1,\ldots,n$ we can apply Theorem \ref{theo difference algebraic Chevalley} to the inclusion $R/\ida_i\subset S/{\ida_i}'$ of $\s^{n'}$-domains to find a non-zero element $\overline{r_i}\in R/\ida_i$ such that every $\s^{dn'}$-prime ideal $\overline{\q_i}$ of $R/\ida_i$ with $\overline{r_i}\notin_{\s^{n'}}\overline{\q_i}$ lifts to a $\s^{l_idn'}$-prime ideal of $S/{\ida_i}'$.

For $i=1,\ldots,n$ let $e_i$ denote an element of $R$ such that
\[e_i\in(\ida_1\cap\cdots\cap\ida_{i-1}\cap\ida_{i+1}\cap\cdots\cap\ida_{n})\smallsetminus\ida_i\]
and set $r=e_1r_1+\cdots+e_nr_n$. As $r_i\notin\ida_i$ we see that $r\notin\ida_i$ for $i=1,\ldots,n$. Thus $r$ is a non-zero divisor of $R$. Let $\p$ be a $\s$-pseudo prime ideal of $R$ with $r\notin_\s\p$. Let $\q$ be a minimal prime ideal of $\p$. Then $\q\supset \ida_i$ for some $i\in\{1,\ldots,n\}$. So $\q$ defines a $\s^{dn'}$-prime ideal $\overline{\q}$ of $R/\ida_i$ for some $d\geq 1$. Suppose $\overline{r_i}\in_{\s^{n'}}\overline{\q}$. Then there exists an integer $j\geq 0$ such $\s^{jn'}(r_i)\in\q$. As $\ida_i\subset \q$ is a $\s^{n'}$-ideal it follows that $\s^{jn'}(e_{i'})\in\q$ for $i'\neq i$ and therefore
\[\s^{jn'}(r)=\s^{jn'}(e_1)\s^{jn'}(r_1)+\cdots+\s^{jn'}(e_i)\s^{jn'}(r_i)+\cdots+\s^{jn'}(e_n)\s^{jn'}(r_n)\]
lies in $\q$. But this contradicts $r\notin_\s\p$.

Thus $\overline{r_i}\notin_{\s^{n'}}\overline{\q}$ and there exists a $\s^{l_idn'}$-prime ideal $\overline{\q'}$ of $S/{\ida_i}'$ such that $\overline{\q'}\cap (R/\ida_i)=\overline{\q}$. Then $\q'$ is a $\s^{l_idn'}$-prime ideal of $S$ with $\q'\cap R=\q$. Consequently the $\s$-pseudo prime ideal $\p'$ of $S$ associated with $\q'$ satisfies $\p'\cap R=\p$. \qed

\vspace{5mm}

It seems natural to ask if the $\s$-algebraic Chevalley theorem presented above will eventually enable one to establish a $\s$-geometric Chevalley theorem, i.e. a difference analog of Theorem \ref{theo usual geometric Chevalley theorem}. However, already the definition of a difference scheme in this context is somewhat beyond the scope of this short note (See however \cite{Hrushovski:elementarytheoryoffrobenius}, \cite{Wibmer:thesis} and \cite{Trushin:DifferenceNullstellsatzCaseOfFiniteGroup}.) A more geometric presentation of the approach indicated in this article remains for the future.

\section{Applications}

\subsection{Constrained extensions}

In differential algebra there is a notion of a differentially closed field and a differential closure of a differential field (see \cite{Kolchin:constrainedExtensions} or the second part in \cite{MarkerMessmerPillay:ModelTheoryofFields}). If $K$ is a differential field
then a finitely differentially generated differential field extension $L$ of $K$ has a $K$-embedding into the differential closure of $K$ if and only if $L$ is a constrained extension of $K$.

In difference algebra there exists a notion of $\s$-algebraically closed $\s$-field (model of ACFA, \cite{Hrushovskietal:ModelTheoryofDifferencefields}) but it seems that there does not exist a satisfactory notion of a $\s$-algebraic closure that is a field. While it is true that every $\s$-field $K$ can be embedded into a model of ACFA, say $\widetilde{K}$, this usually has some shortcomings: First, due to the phenomenon of incompatibility, a system of $\s$-algebraic equations with coefficients in $K$ that has a solution in a $\s$-overfield of $K$ need not have a solution in $\widetilde{K}$. Secondly, if the constants of $K$ are algebraically closed then $\widetilde{K}$ will have new constants, whereas in the differential case the differential closure of a differential field with algebraically closed constants will not have new constants.

In this section we propose a difference analog of the constrained extensions in differential algebra and establish their basic properties. Philosophically finitely $\s$-generated constrained extensions should play the same role in difference algebraic geometry as finite algebraic extensions in algebraic geometry.


\begin{defi}
A $\s$-pseudo domain is called \emph{$\s$-pseudo simple} if it has no $\s$-pseudo prime ideal other than the zero ideal.
\end{defi}

In some sense the $\s$-geometrical interpretation of a $\s$-pseudo simple ring is a point.
\begin{rem} \label{rem pseudo simple} Let $R$ be a $\s$-pseudo domain.
\begin{enumerate}
\item If $R$ is noetherian then $R$ is $\s$-pseudo simple if and only if $R$ is $\s$-simple.
\item $R$ is $\s$-pseudo simple if and only if $R$ has no proper $\s^d$-perfect ideal (for all $d\geq 1$).
\end{enumerate}
\end{rem}
\noindent Proof: (i): Let $R$ be a noetherian $\s$-pseudo simple $\s$-domain and $\ida$ a $\s$-ideal of $R$. The radical of the reflexive closure of $\ida$ is a proper reflexive radical $\s$-ideal of $R$. As $\ida$ is reflexive we see that $\q\mapsto\s^{-1}(\q)$ induces a bijection on the finite set of minimal prime ideals of $\ida$. Therefore $\ida$ is the intersection of $\s$-pseudo prime ideals. Thus $\ida=0$.

(ii): This is immediate from the fact that a $\s$-perfect ideal is the intersection of $\s$-prime ideals (See \cite[p. 88]{Cohn:difference}). \qed

\begin{defi}
Let $L|K$ be an extension of $\s$-pseudo fields. A tuple $a\in L^n$ is called\emph{ constrained over $K$} if there exists a non-zero divisor $b\in K\{a\}$ such that $K\{a,\frac{1}{b}\}$ is $\s$-pseudo simple. In this situation we also call $b$ a \emph{constraint} of $a$ over $K$.
The extension $L$ is called constrained over $K$ if every finite tuple of elements in $L$ is constrained over $K$.
\end{defi}

If $R$ is a $\s$-pseudo domain containing a $\s$-pseudo field $K$, we also say that $R$ is constrained over $K$ if $\Q(R)$ is constrained over $K$.

\begin{ex}
If $L|K$ is an extension of $\s$-fields such that $L$ is algebraic over $K$ then $L$ is constrained over $K$.
\end{ex}
\begin{ex}\label{ex Picard-Vessiot}
Let $K$ be an inversive $\s$-field with algebraically closed constants. Then a total Picard-Vessiot ring over $K$ (in the sense of \cite[Definition 1.22, p. 16]{SingerPut:difference}) is a constrained $\s$-pseudo field extension of $K$.
\end{ex}
\noindent Proof: By definition a total Picard-Vessiot ring is the total quotient ring of a finitely generated, $\s$-simple $K$-$\s$-algebra. Thus the claim follows from Proposition \ref{prop constrained welldefined} below.\qed

\begin{ex}
A generalization of Example \ref{ex Picard-Vessiot} is the fact that a $\s$-pseudo field extension satisfying the conditions analogous to strongly normal differential field extensions is constrained. See \cite[Lemma 3.4.1, p. 79]{Wibmer:thesis}.
\end{ex}

\begin{ex}\label{ex constrained zykel}
Let $K$ be a $\s$-field and $d\geq 1$ an integer. Set $L=K\oplus\cdots\oplus K$ and define $\s: L\rightarrow L$ by
\[\s(a_1,\ldots,a_d)=(\s(a_d),\s(a_1),\ldots,\s(a_{d-1})).\]
Then $L$ is a constrained $\s$-pseudo field extension of $K$. ($K$ is diagonally embedded into $L$.) A $\s$-pseudo field extension of this form will be called \emph{trivial}.
\end{ex}

\begin{ex} \label{ex constrained extensions of ACFA}
Let $K$ be a model of ACFA. Then every finitely $\s$-generated constrained $\s$-pseudo field extension of $K$ is trivial.
\end{ex}
\noindent Proof: Let $L$ be a finitely $\s$-generated constrained $\s$-pseudo field extension of $K$. Then $L=k(\p)$ for some $\s$-pseudo prime ideal $\p$ of some $\s$-polynomial ring $K\{x\}$ such that there is no $\s$-pseudo prime ideal of $K\{x\}$ properly containing $\p$. Let $\q$ be a minimal prime ideal of $\p$ and let $d$ denote the period of $\q$. We can interpret $\q$ as a $\s^d$-prime ideal in the $\s^d$-polynomial ring $K\{x\}=K\{x,\s(x),\ldots,\s^{d-1}(x)\}_{\s^d}$ over the $\s^d$-field $K$.
By \cite[Corollary 1.12, p. 3013]{Hrushovskietal:ModelTheoryofDifferencefields} we know that $(K,\s^d)$ also is a model of ACFA. As there is no $\s^d$-prime ideal strictly above $\q$ this shows that $k(\q)=K$. \qed

\begin{prop}\label{prop constrained welldefined}
Let $L|K$ be an extension of $\s$-pseudo fields and $a\in L^n$. If $a$ is constrained over $K$ then $K\langle a\rangle$ is constrained over $K$.
\end{prop}
\noindent Proof: Let $a'\in K\langle a\rangle^{n'}$. There exists a non-zero divisor $c\in K\{a\}$ such that $K\{a'\}\subset K\{a,\frac{1}{c}\}$. As $a$ is constrained over $K$ there exists a non-zero divisor $b\in K\{a\}$ such that $K\{a,\frac{1}{b}\}$ is $\s$-pseudo simple. Then also $K\{a,\frac{1}{b},\frac{1}{c}\}$ is $\s$-pseudo simple. If we apply Corollary \ref{cor chevalley for pseudo domain} to the inclusion $K\{a'\}\subset K\{a,\frac{1}{b},\frac{1}{c}\}$ we find that there exists a non-zero divisor $r\in K\{a'\}$ such that every $\s$-pseudo prime ideal $\p$ of $K\{a'\}$ with $r\notin_\s\p$ lifts to $K\{a,\frac{1}{b},\frac{1}{c}\}$. But the zero ideal is the only $\s$-pseudo prime ideal of $K\{a,\frac{1}{b},\frac{1}{c}\}$. Thus the only $\s$-pseudo prime ideal of $K\{a'\}$ with $r\notin_\s\p$ is the zero ideal. Hence $R\{a',\frac{1}{r}\}$ is $\s$-pseudo simple, i.e., $a'$ is constrained with constraint $r$. \qed

\begin{prop} \label{prop constrained compatibel}
Let $L|K$ and $M|L$ be finitely $\s$-generated extensions of $\s$-pseudo fields.
\begin{enumerate}
\item If $L$ is constrained over $K$ and $M$ is constrained over $L$ then $M$ is constrained over $K$.
\item If $M$ is constrained over $K$ then $L$ is constrained over $K$ and $M$ is constrained over $L$.
\end{enumerate}
\end{prop}
\noindent Proof: (i): Let $a\in M^n$. We have to show that $a$ is constrained over $K$. As $a$ is constrained over $L$ there exists a non-zero divisor $b\in L\{a\}$ such that $L\{a,\frac{1}{b}\}$ is $\s$-pseudo simple. There is a finite $L$-tuple $c$ such that $K\{c\}$ is $\s$-pseudo simple and $K\langle c\rangle =L$.
Then $K\{c,a,\frac{1}{b}\}$ is $\s$-pseudo simple (as it has the same pseudo spectrum as $L\{a,\frac{1}{b}\}$). Applying Corollary \ref{cor chevalley for pseudo domain} to the
inclusion $K\{a\}\subset K\{c,a,\frac{1}{b}\}$ shows that $a$ is constrained over $K$.

(ii): We assume that $M$ is constrained over $K$. It is trivial that $L$ is constrained over $K$, so we only have to show that $M$ is constrained over $L$. Let $a\in M^n$. Let $b$ be finite $\s$-generating set of $L$ over $K$. As the $M$-tuple $(b,a)$ is constrained over $K$ there exists a non-zero divisor $c\in K\{b,a\}$ such that $K\{b,a,\frac{1}{c}\}$ is $\s$-pseudo simple. But then $c\in L\{a\}$ and also $L\{a,\frac{1}{c}\}$ is $\s$-pseudo simple, i.e., $a$ is constrained over $L$. \qed

\begin{prop}
Let $L|K$ be a constrained extension of $\s$-pseudo fields. Then the constants of $L$ are algebraic over the constants of $K$. In particular, if the constants of $K$ are algebraically closed then $L$ does not have new constants.
\end{prop}
\noindent Proof: With the aid of \cite[Lemma 1.3.5, p. 10]{Wibmer:thesis} we can easily reduce to the case that $K$ is a field. Let $C=K^\s$ denote the constants of $K$ and let $c\in L$ be constant.
Suppose that $c$ is transcendental over $C$. As $c$ is constrained over $K$ there exists a non-zero divisor $b\in K\{c\}$ such that $K\{c,\frac{1}{b}\}$ is $\s$-pseudo simple. In other words, every $\s$-pseudo prime ideal $\p$ of $K\{c\}$ with $b\notin_\s\p$ is zero.
One knows that $K$ and $L^\s$ are linearly disjoint over $C$ (see e.g. \cite[Lemma 1.1.6, p. 4]{Wibmer:thesis}). Hence $K\{c\}=K[c]=K\otimes_C C[c]$ is a univariate polynomial ring over the field $K$.
From \cite[Proposition 1.4.15. p. 15]{Wibmer:thesis} we know that the $\s$-pseudo prime ideals of $K\{c\}$ are precisely those of the form $K\otimes \q$ where $\q$ is a prime ideal of $C[c]$. The prime ideals of $C[c]$ are given by irreducible polynomials. If $f\in C[c]$ is such an irreducible polynomial then $f$ factors as a finite product of
irreducible polynomials in $K[c]$. We have $b\in_\s K\otimes (f)\subset K\{c\}$ if $b$ is divisible by one of these factors. Distinct (monic) irreducible polynomials of $C[c]$ do not give rise to common factors. Thus, for all but finitely many irreducible polynomials $(f)$ of $C[c]$ we have $b\notin_\s K\otimes (f)$. But then $K\otimes (f)$ must be zero -- a contradiction. \qed

%
%
%

\vspace{5mm}

The following proposition (together with Proposition \ref{prop constrained welldefined}) provides a ``generic'' way of constructing constrained extensions.

\begin{prop} \label{prop exists maximal pseudo prime}
Let $K$ be a $\s$-pseudo field and $R$ a $K$-$\s$-algebra that is finitely $\s$-generated over $K$. Assume that $R$ has a $\s$-pseudo prime ideal. Then there exists a maximal element in the set of all $\s$-pseudo prime ideals of $R$ ordered by inclusion.
\end{prop}
\noindent Proof: We can easily reduce to the case that $K$ is a field.
We have to show that there is a maximal element in the set of all $\s^\bullet$-prime ideals of $K\{x\}=K\{x_1,\ldots,x_n\}$ containing a given $\s^d$-prime ideal $\q$ of $K\{x\}$.
The obvious application of Zorn's Lemma shows that there exists a maximal element, say $\q'$, in the set of all $\s^d$-prime ideals of $K\{x\}$ that contain $\q$.

We want to show that the transcendence degree of $\Q(K\{x\}/\q')$ over $K$ is finite. We can interpret $\q'$ as a $\s^d$-prime ideal in the $\s^d$-polynomial ring $K\{x\}_\s=K\{x,\ldots,\s^{d-1}(x)\}_{\s^d}$. By construction the $\s^d$-variety defined by $\q'$ has only one point. (That is, every solution of $\q'$ in a $\s^d$-field extension of $K$ is equivalent to the generic solution.) Because every difference algebraic variety has a point in a difference algebraic extensions of the ground field
(see \cite[Chapter 8, Section 2]{Cohn:difference}) it follows that $\Q(K\{x\}/\q')$ is a $\s^d$-algebraic extension of $K$. In particular, the transcendence degree of $\Q(K\{x\}/\q')$ is finite.

In general the Krull dimension of a $K$-domain is at most the transcendence degree of its quotient field over $K$ (\cite{Giral:KrullDimension}). Therefore any strictly increasing chain of $\s^\bullet$-prime ideals above $\q'$ must be finite. \qed


\vspace{5mm}

The following corollary is the difference version of \cite[Proposition 6, p. 142]{Kolchin:differentialalgebraandalgebraicgroups}. But please note that the proof for the difference case is much more complicated.

\begin{cor}
Let $L|K$ be an extension of $\s$-pseudo fields. Let $a\in L^n$ and $b\in K\{ a\}$ a non-zero divisor. Then there exists a $K$-morphism $\psi:K\{a\}\rightarrow L'$ into a
$\s$-pseudo field extension $L'$ of $K$ such that $\psi(a)$ is constrained over $K$ with constraint $\psi(b)$.
\end{cor}
\noindent Proof: By Proposition \ref{prop exists maximal pseudo prime} there exists a maximal element, say $\p$, in the set of all $\s$-pseudo prime ideals of $K\{a,\frac{1}{b}\}$. Set
$L'=\Q\left(K\{a,\frac{1}{b}\}/\p\right)$ and let $\psi:K\{a\}\rightarrow L$ denote the canonical map. Then $L'$ and $\psi$ have the desired properties. \qed

\vspace{5mm}

Maybe the following reformulation is somewhat more catchy:
\begin{theo} \label{theo have solution constrained}
If a system of algebraic difference equations with coefficients in a $\s$-pseudo field $K$ has a solution in a $\s$-pseudo field extension of $K$ then it already has a solution in a finitely $\s$-generated constrained $\s$-pseudo field extension of $K$. \qed
\end{theo}

The above theorem confirms the point of view that finitely $\s$-generated constrained extensions of a $\s$-field are the analog of finite algebraic extensions of the ground field in usual algebraic geometry.

\subsection{Differential Picard-Vessiot theory with a difference parameter}

A differential Galois theory with differential parameters has been developed in the general strongly normal context in \cite{Landesman:GeneralizedDifferentialGaloisTheory}. The special case of linear differential equations with differential parameters is presented in more detail in \cite{CassidySinger:GaloisTheoryofParameterizedDifferentialEquations}. The approach of \cite{CassidySinger:GaloisTheoryofParameterizedDifferentialEquations} has been generalized to include linear difference equations with differential parameters in \cite{HardouinSinger:DifferentialGaloisTheoryofLinearDifferenceEquations}.

Recently a Galois theory for linear difference equations with difference parameters has been proposed in \cite{OvchinnikovTrushinetal:GaloisTheory_of_difference_equations_with_periodic_parameters}. However the approach of \cite{OvchinnikovTrushinetal:GaloisTheory_of_difference_equations_with_periodic_parameters} only applies to parameters of finite order.

Currently a Galois theory for linear differential equations with difference parameters is being developed in \cite{HardouindiVizio:DifferenceGaloisofDifferential}. Here we confine to present the connection to constrained extensions.

In the usual Picard-Vessiot theory one can show that the Picard-Vessiot ring of a linear differential equation is unique up to a finite algebraic extension of the differential constants. To support our claim that finitely $\s$-generated constrained extensions are the difference counterpart of finite algebraic extensions we shall prove that the Picard-Vessiot ring of a linear differential equation with a difference parameter is unique up to a finitely $\s$-generated constrained extension of the differential constants.

\vspace{5mm}

First we need to fix some notation. By a $\de\s$-ring $R$ we mean a ring equipped with a derivation $\de:R\rightarrow R$ and an endomorphism $\s:R\rightarrow R$ such that $\de$ and $\s$ commute. This implies that the $\de$-constants $R^\de=\{r\in R| \ \de(r)=0\}$ are a $\s$-ring.

Throughout this section we denote with $k$ a $\de\s$-field. For simplicity we assume that $k$ is of characteristic zero \footnote{In general one should use iterative derivations.}.
We denote with $C=k^\de$ the $\s$-field of $\de$-constants of $k$. We also assume that $k$ is $\s$-separable over $C$ (see \cite[Section 1.5, p. 16]{Wibmer:thesis}). This means that if a subset $S$ of $k$ is linearly independent over $C$ then also $\s(S)$ is linearly independent over $C$. This is always satisfied if $\s:C\rightarrow C$ is surjective, which in turn is the case if $\s:k\rightarrow k$ is surjective.

\vspace{5mm}

\begin{defi}
Let $A\in k^{n\times n}$ and $R$ a $k$-$\delta\sigma$-algebra with the property that there exists $Y\in \Gl_n(R)$ such that $\de(Y)=AY$ and $R=k\left\{Y_{ij},\frac{1}{\det(Y)}\right\}_\s$. Then $R$ is called a \emph{Picard-Vessiot ring (with $\s$-parameter) for $A$} if $R$ is $\de$-simple.
\end{defi}

We stress the point that we require a Picard-Vessiot ring to be $\de$-simple and not just $\de\s$-simple.
This has the advantage that a Picard-Vessiot ring is automatically an integral domain because a $\de$-simple ring is an integral domain.
By \cite[Corollary 6.22, p. 372]{HardouinSinger:DifferentialGaloisTheoryofLinearDifferenceEquations} the requirement that $R$ is  $\de$-simple is conform with the standard approach. Since the $\s$-radical of $R$, i.e., $\{r\in R|\ \exists \  n\geq 1\text{ such that } \s^n(r)=0\}$ is a $\de$-ideal it also follows that $R$ is a $\s$-domain.
Moreover a Picard-Vessiot ring is $\de\s$-simple.

As the first step one needs the existence of Picard-Vessiot rings. It is important to notice that the usual strategy for constructing a Picard-Vessiot ring fails in this context. That is, taking the quotient by a $\de\s$-maximal ideal in the generic solution ring $k\{X,\frac{1}{\det(X)}\}_\s$ does not necessarily yield something useful. It is not even clear to the author if such an ideal will have only finitely many minimal prime ideals. However, a more sensitive adaption of the usual construction to the difference world works well enough.

\begin{lemma}
Let $A\in k^{n\times n}$. Then there exists a Picard-Vessiot ring with $\s$-parameter for $A$.
\end{lemma}
\noindent Proof: The proof is analogous to the construction of prolongations of difference kernels (see \cite[Chapter 6]{Cohn:difference} or Section 1.3 above). Let $X$ be a $n\times n$-tuple of $\s$-indeterminates over $k$. We consider the generic solution ring $S=k\{X,\frac{1}{\det(X)}\}_\s$. To be precise: $S$ is the localization of the $\s$-polynomial ring in the $X_{ij}$'s ($1\leq i,j\leq n$) over $k$ at the multiplicatively closed subset generated by $\det(X),\s(\det(X)),\ldots$. We consider $S$ as $\de$-ring by setting $\de(X)=AX, \de(\s(X))=\s(A)\s(X),\ldots$. Then $S$ is a $k$-$\de\s$-algebra containing the fundamental solution matrix $X$ for $A$. For $d\geq 0$ we abbreviate
\[S_d=k\left[X,\s(X),\ldots,\s^d(X),\frac{1}{\det(X)},\ldots,\s^d\left(\frac{1}{\det(X)}\right)\right]\subset S.\]
Then $S_d$ is a $k$-$\de$-algebra and $\s$ restricts to $\s: S_{d-1}\rightarrow S_{d}$. We will prove by induction on $d\geq 0$ that there exists a $\de$-maximal $\de$-ideal $\q_d$ of $S_d$ such that the inverse image of $\q_d$ under $\s: S_{d-1}\rightarrow S_{d}$ equals $\q_{d-1}=\q_d\cap S_{d-1}$.
For $d=0$ we simply choose a $\de$-maximal ideal $\q_0$ of $S_0$. (Note that $S_0/\q_0$ is a usual Picard-Vessiot ring for $A$ over the $\de$-field $k$.)  Now we do the induction step: Let $\ida$ denote the ideal of $S_d$ generated by $\q_{d-1}$ and $\s(\q_{d-1})$. Then $\ida$ is a $\de$-ideal. By the theory of difference kernels $\q_{d-1}$ can be prolonged, in particular $1\notin\ida$. Thus there exists a $\de$-maximal $\de$-ideal $\q_d$ of $S_d$ containing $\ida$. The inverse image of $\q_d$ under $\s:S_{d-1}\rightarrow S_d$ is a $\de$-ideal containing $\q_{d-1}$. As $\q_{d-1}$ is $\de$-maximal it follows that the inverse image of $\q_d$ under $\s:S_{d-1}\rightarrow S_d$ equals $\q_{d-1}$.
Similarly $S_{d-1}\cap\q_d$ is a $\de$-ideal of $S_{d-1}$ containing $\q_{d-1}$. Therefore, the inverse image of $\q_d$ under $\s:S_{d-1}\rightarrow S_d$ equals $\q_{d-1}=\q_d\cap S_{d-1}$.

We conclude that $\q=\cup\q_d$ is a $\s$-prime $\de\s$-ideal of $S$. Clearly $\q$ is $\de$-maximal: If $\m$ was a proper $\de$-ideal of $S$ properly containing $\q$ then, for some $d$, $\m\cap S_d$ would be a $\de$-ideal properly containing $\q_d=\q\cap S_d$ -- in contradiction to the $\de$-maximality of $\q_d$.
Therefore $R=S/\q$ is a Picard-Vessiot ring for $A$. \qed

%
%
%


\vspace{5mm}

The key tool for our task is the following proposition.

\begin{prop} \label{prop constants constrained}
Let $R$ be a $k$-$\de\s$-algebra such that $R$ is $\de$-simple and finitely $\s$-generated over $k$. Then $R^\de$ is a finitely $\s$-generated constrained $\s$-field extension of $C=k^\de$.
\end{prop}
\noindent
Proof: First of all we note that $R^\de$ is a field because $R$ is $\de$-simple. Let $c$ be a finite tuple with coordinates in $R^\de$. If we apply Theorem \ref{theo difference algebraic Chevalley} to the inclusion
$k\{c\}_{\s}\subset R$ we find that there exists an integer $l\geq 1$ and a non-zero $r$ in $k\{c\}_{\s}$ such that every $\s^d$-prime ideal $\widetilde{\q}$ of $k\{c\}_{\s}$ with $r\notin_\s\widetilde{\q}$ lifts to a $\s^{ld}$-prime ideal of $R$.

We may write $r=\lambda_1\otimes a_1+\cdots +\lambda_n\otimes a_n\in k\otimes_C C\{c\}_\s=k\{c\}_\s$ with the $\lambda_i$'s linearly independent over $C$.
Let $b\in C\{c\}_\s$ denote a non-zero $a_i$. We claim that $b$ is a constraint for $c$ over $C$.
Let $\q$ be a non-zero $\s^d$-prime ideal of $C\{c\}_\s$ with $b\notin_\s\q$. We have to show that $\q=0$.

As $C=k^\de$ is relatively algebraically closed in $k$ (\cite[Corollary, p. 94]{Kolchin:differentialalgebraandalgebraicgroups}) one knows
(e.g. \cite[Proposition 9, A.V.142]{Bourbaki:Algebra2}) that $k\otimes\q$ is a prime ideal of $k\otimes_C C\{c\}_\s$.
By assumption $k$ is $\s$-separable and also $\s^d$-separable over $C$. It follows from \cite[Corollary 1.5.4, p. 18]{Wibmer:thesis} that $k\otimes\q$ is a $\s^d$-prime ideal of $k\{c\}_\s$.
Suppose, for a contradiction, that $r\in_\s k\otimes\q$. Then exists $m\geq 0$ such that $\s^m(r)\in k\otimes\q$. Because $k$ is $\s$-separable over $C$ the $\s^m(\lambda_i)$ are $C$-linearly independent. By considering the image of $\s^m(r)$ in
\[\left(k\otimes_C C\{c\}_\s\right)/k\otimes\q=k\otimes_C\left( C\{c\}_\s/\q\right)\]
we see that this implies $\s^{m}(b)\in\q$ -- in contradiction to $b\notin_\s\q$.

Hence $r\notin\widetilde{\q}=k\otimes \q$ and so $\widetilde{\q}$ lifts to a $\s^{ld}$-prime ideal $\q'$ of $R$. In particular $\q'$ contains the ideal generated by $\q$ in $R$, but because $\q\subset R^\de$ this ideal is a $\de$-ideal. Thus it must be zero, and so $\q=0$ as desired.

\vspace{5mm}

It remains to see that $R^\de$ is finitely $\s$-generated over $C$. As $\Q(R)$ is finitely $\s$-generated over $k$ it follows that $\Q(k\otimes_C R^\de)$ is finitely $\s$-generated over $k$ (\cite[Theorem 18, p. 145]{Cohn:difference}). There exists a finite $\s$-generating set $S$ of $\Q(k\otimes_C R^\de)$ as $\s$-field extension of $k$ consisting of elements of $R^\de$. Then $S$ is also a $\s$-generating set of $R^\de$ as $\s$-field extension of $C$. \qed

\begin{cor} \label{cor no new constants}
Let $R$ be a Picard-Vessiot ring. Then $R^\de$ is a finitely $\s$-generated constrained $\s$-field extension of $C=k^\de$. In particular if $C$ does not have proper constrained $\s$-field extensions (e.g. $C$ is a Model of ACFA) then $R$ does not have new $\de$-constants. \qed
\end{cor}

Now we are prepared to prove the general uniqueness theorem for differential Picard-Vessiot rings with a difference parameter.
\begin{theo} \label{theo uniqueness}
Let $R_1$ and $R_2$ be two Picard-Vessiot rings for $A\in k^{n\times n}$. Then there exists a finitely $\s$-generated constrained $\s$-pseudo field extension $\widetilde{C}$ of $C$ containing $C_1=R_1^\de$ and $C_2=R_2^\de$ such that $R_1\otimes_{C_1} \widetilde{C}$ and $R_2\otimes_{C_2}\widetilde{C}$ are isomorphic as $k\otimes_C \widetilde{C}$-$\de\s$-algebras.
\end{theo}
\noindent Proof: We consider the $\de\s$-$k$-algebra $R_1\otimes_k R_2$. Let $Y_1\in \Gl_n(R_1)$ and $Y_2\in\Gl_n(R_2)$ denote fundamental matrices for $A$.
Let $D=(Y_1\otimes 1)^{-1}(1 \otimes Y_2)\in\Gl_n(R_1\otimes_k R_2)$. A well known computation shows that $\de(D)=0$.
Because $1\otimes Y_2=Y_1\otimes 1 D$ and $R_2$ is $\s$-generated by $Y_2$ and $\frac{1}{\det(Y_2)}$ it follows that
\[R_1\otimes_k R_2=R_1\otimes_{C_1}C_1\Big\{D,\frac{1}{\det(D)}\Big\}_\s.\]

Let $L_1=\Q(R_1)$ and $L_2=\Q(R_2)$. Then $L_1$ and $L_2$ are finitely $\s$-generated $\s$-field extensions of $k$. It follows from Theorem \ref{theo compatibility} that there exists a $\s$-pseudo prime ideal in $L_1\otimes_k L_2$. Thus there also is a $\s$-pseudo prime ideal in $R_1\otimes_kR_2$ and $C_1\{D,\frac{1}{\det(D)}\}_\s$.
Hence by Theorem \ref{theo have solution constrained} there exists a $C_1$-$\s$-morphism $\psi:C_1\{D,\frac{1}{\det(D)}\}_\s\rightarrow \widetilde{C}$ into a finitely $\s$-generated constrained $\s$-pseudo field extension $\widetilde{C}$ of $C_1$. From Corollary \ref{cor no new constants} we know that $C_1$ is a finitely $\s$-generated constrained $\s$-pseudo field extension of $C$. Thus $\widetilde{C}$ is finitely $\s$-generated over $C$ and from Proposition \ref{prop constrained compatibel} (i) it follows that $\widetilde{C}$ is constrained over $C$.

We have a $k$-$\de\s$-morphism
\[\phi:R_2\longrightarrow R_1\otimes_k R_2=R_1\otimes_{C_1}C_1\Big\{D,\frac{1}{\det(D)}\Big\}_\s\xrightarrow{\id\otimes \psi} R_1\otimes_{C_1}\widetilde{C}.\]
As $C_2=R_2^\de\subset(R_1\otimes_k R_2)^\de=C_1\{D,\frac{1}{\det(D)}\}_\s$ we have a natural inclusion of $C_2$ in $\widetilde{C}$ and we can extend $\phi$ to a $k\otimes_C \widetilde{C}$-$\de\s$ morphism
\[\phi:R_2\otimes_{C_2}\widetilde{C}\longrightarrow R_1\otimes_{C_1}\widetilde{C}.\]
We will show that $\phi$ is an isomorphism. For the surjectivity it suffices to show that the entries of $Y_1$ and $\frac{1}{\det(Y_1)}$ lie in the image of $\phi$. But this is clear from $Y_1\otimes 1=(1\otimes Y_2) D^{-1}$. To show that $\phi$ is injective it suffices to notice that $R_2\otimes_{C_2}\widetilde{C}$ is $\de\s$-simple by Lemma \ref{lemma deltasigma simple} below.
\qed

\begin{lemma}  \label{lemma deltasigma simple}
Let $R$ be a $\de$-simple $k$-$\de\s$-algebra and $\widetilde{C}$ a $\s$-pseudo field extension of $C'=R^\de$ (considered as constant $\de$-ring). Then
$R\otimes_{C'}\widetilde{C}$ is $\de\s$-simple.
\end{lemma}
\noindent Proof: Every $\de$-ideal of $R\otimes_{C'}\widetilde{C}$ is of the form $R\otimes \ida$ for some ideal $\ida$ of $\widetilde{C}$. This follows
e.g. as in the proof of \cite[Proposition 5.6, p. 4484]{Kovacic:differentialgaloistheoryofstronglynormal}. But if $\ida$ is a non-zero ideal of $\widetilde{C}$ then it contains a minimal idempotent of $\widetilde{C}$. Therefore if $R\otimes \ida$ is a $\s$-ideal it must contain $1$. \qed

\vspace{5mm}

We recall that the definition of a trivial $\s$-pseudo field extension was given in Example \ref{ex constrained zykel}. And in Example \ref{ex constrained extensions of ACFA} we noted that every finitely $\s$-generated constrained $\s$-pseudo field extension of a model of ACFA is trivial.

\begin{cor} \label{cor uniqueness}
Let $k$ be a $\de\s$-field and assume that all finitely $\s$-generated constrained $\s$-pseudo field extensions of $k^\de$ are trivial (e.g. $k^\de$ is a model of ACFA). Let $R_1$ and $R_2$ be Picard-Vessiot rings with a $\s$-parameter for the differential equation $\de(Y)=AY$ with $A\in k^{n\times n}$. Then there exists an integer $l\geq 1$ such that $R_1$ and $R_2$ are isomorphic as $k$-$\de\s^l$-algebras.
\end{cor}
\noindent Proof: By Corollary \ref{cor no new constants} we have $R_1^\de=C$ and $R_2^\de=C$ with $C=k^\de$. By Theorem \ref{theo uniqueness} and the assumptions it follows that there is an integer $l\geq 1$ such that $R_1\oplus\cdots\oplus R_1$ and $R_2\oplus\cdots\oplus R_2$ are $k$-$\de\s$-isomorphic. Here the sums have $l$ summands, the derivation is given by
$\de(a_1\oplus\cdots\oplus a_l)=\de(a_1)\oplus\cdots\oplus\de(a_l)$ and $\s$ is given by $\s(a_1\oplus\cdots\oplus a_l)=\s(a_l)\oplus\s(a_1)\oplus\cdots\oplus\s(a_{l-1})$.\qed

\vspace{5mm}

The following simple example shows that one can not in general choose $l=1$ in the above corollary.
\begin{ex}
Let $C$ be a model of ACFA of characteristic zero. Let $k=C(x)$ denote the rational functions in a variable $x$ over $C$. Choose $q\in C\smallsetminus\{0\}$. We consider $k$ as $\de\s$-field by setting $\de=x\frac{d}{dx}$ and $\s(x)=qx$. Then $k^\de=C$. Consider the equation $\de(y)=\frac{1}{2}y$. The derivation $\de$ uniquely extends to
$C(\sqrt{x})\supset C(x)=k$ by $\de(\sqrt{x})=\frac{1}{2}\sqrt{x}$, so that $\sqrt{x}$ is a fundamental solution matrix for $\de(y)=\frac{1}{2}y$. For the extension of $\s$ we have two choices $\s(\sqrt{x})=\sqrt{q}\sqrt{x}$ or $\s(\sqrt{x})=-\sqrt{q}\sqrt{x}$. Both choices turn $C(\sqrt{x})$ into a Picard-Vessiot ring with $\s$-parameter but clearly the difference structures are not isomorphic, however, as predicted by Corollary \ref{cor uniqueness} they become isomorphic if we pass from $\s$ to $\s^2$.
\end{ex}

Please note that the assumption that $C$ is a model of ACFA is completely irrelevant to the above example. This means that whatever assumption one is willing to make on the $\de$-constants of $k$, it will not suffice to guarantee uniqueness in the stringent sense of the word in general.

Nevertheless we have the following uniqueness result result. (Recall that a $\s$-field $k$ is called universally compatible if any two $\s$-field extensions of $k$ are compatible.)

\begin{cor}
Let $k$ be a $\de\s$-field and assume that $k^\de$ is a model of ACFA. If $k$ is universally compatible as $\s$-field (e.g. $k$ is algebraically closed), then a Picard-Vessiot ring with $\s$-parameter is unique (up to $k$-$\de\s$-isomorphisms.)
\end{cor}
\noindent Proof: We use the notation of the proof of Theorem \ref{theo uniqueness}. By Corollary \ref{cor no new constants} we have $C_1=C_2=C$. By assumption $L_1|k$ and $L_2|k$ are compatible. This implies that there is a $\s$-prime ideal in $C\{D,\frac{1}{\det(D)}\}_\s$. By the definition of ACFA there exists a $C$-$\s$-morphism $\psi:C\{D,\frac{1}{\det(D)}\}_\s\rightarrow C$. In other words, we can choose $\widetilde{C}=C$. \qed

\vspace{5mm}

Finally we note that the twist (i.e. the passage from $\s$ to $\s^l$) required for the uniqueness in Corollary \ref{cor uniqueness} poses no serious problem for defining the Galois group of the equation, which, after all is the main purpose of the Picard-Vessiot ring. If $R_1$ and $R_2$ are $\de\s^l$-isomorphic then we have an induced $k$-$\de\s^l$-isomorphism $R_1\otimes_k R_1\rightarrow R_2\otimes_k R_2$ which induces a $C$-$\s^l$ isomorphism $(R_1\otimes_k R_1)^\de\rightarrow (R_2\otimes_k R_2)^\de$ between the $\s$-coordinate rings of the corresponding Galois groups (which are of course $\s$-algebraic groups).
The $\s^\bullet$-spectrum of a difference ring equals the $\s^{l\bullet}$-spectrum and so $(R_1\otimes_k R_1)^\de$ and $(R_2\otimes_k R_2)^\de$ have the same $\s^\bullet$-spectrum (at least set-theoretically).

%

\bibliographystyle{plain}
\bibliography{bibdata.bib}

\end{document}